\documentclass[12pt]{article}
\usepackage{amsmath,amssymb,amsthm}

\title{ Orbit equivalence for Cantor minimal $\Z^{2}$-systems}
\author{Thierry Giordano\thanks{Supported in part by a
grant from NSERC, Canada}, \\
Department of Mathematics and Statistics, \\
University of Ottawa,\\
585 King Edward, Ottawa, Ontario, Canada K1N 6N5
\and
Hiroki Matui\thanks{Supported in part by a 
grant from the Japan Society for the Promotion of Science},\\
Graduate School of Science and Technology,\\
Chiba University, \\
1-33 Yayoi-cho, Inage-ku, \\
Chiba 263-8522, Japan
 \and
 Ian F. Putnam\thanks{Supported in part by a
grant from NSERC, Canada},\\
Department of Mathematics and Statistics,\\
University of Victoria,\\
Victoria, B.C., Canada V8W 3P4
\and 
Christian F. Skau\thanks{Supported in
part by the Norwegian Research Council}, \\
Department of Mathematical Sciences, \\
Norwegian
University of Science and Technology (NTNU), \\
N-7491 Trondheim, Norway}

\date{  }

\newcommand{\Z}{\mathbb{Z}}

\newcommand{\R}{\mathbb{R}}

\newtheorem{defn}{Definition}[section]
\newtheorem{thm}[defn]{Theorem}

\newtheorem{prop}[defn]{Proposition}
\newtheorem{cor}[defn]{Corollary}
\newtheorem{lemma}[defn]{Lemma}

\begin{document}

\maketitle

\pagebreak

\begin{abstract}
We show that every minimal, free action of the group $\Z^{2}$ on the Cantor set
is orbit equivalent to an AF-relation. As a consequence, this extends the 
classification of minimal systems on the Cantor set up to orbit
equivalence to include AF-relations, $\Z$-actions and $\Z^{2}$-actions.
\end{abstract}

\section{Introduction and statement of results}

In this paper, we consider dynamical systems on the Cantor set. By a Cantor set, $X$, we mean a metrizable
topological space which is compact, totally disconnected (the closed and open sets form a base for the topology) 
and has no isolated points. Any two such spaces are homeomorphic. Moreover, dynamical systems 
on such spaces, which include many symbolic systems, have a fundamental r\^{o}le in the theory.

Our results will be concerned with actions of the group $\Z^{2}$ on
a Cantor set $X$. For the moment, it will be convenient to
consider the case that $G$ is a countable abelian 
group with an action, $\varphi$ on $X$ by homeomorphisms.
For each $g$ in $G$, we have a homeomorphism $\varphi^{g}: X \rightarrow X$ such that, 
$\varphi^{0}(x) =x$, for all $x$ in $X$ and $\varphi^{g+h} = \varphi^{g} \circ \varphi^{h}$.
Furthermore, we will assume the action is free, meaning $\varphi^{g}(x) = x$, for some $x$ in $X$,
 only if $g=0$, 
and minimal, meaning that, for every $x$ in $X$, the set $\{ \varphi^{g}(x) \mid g \in G \}$ is dense
in $X$, or equivalently, there are no non-trivial closed $\varphi$-invariant subsets of $X$.

In addition to group actions, we also consider equivalence relations $R$ on
$X$ which are equipped with a topology in which they are \'{e}tale. The reader should see \cite{Ren:book,PPZ:survey,GPS:affable} for
more information. Roughly speaking, the key idea is that, in this topology, $R$ is locally compact 
and Hausdorff and the two canonical projections
 from $R$ to $X$ are open and, locally, are homeomorphisms. It is important
to realize that this topology is rarely the relative topology from $R \subset X \times X$, except
in the special case that $R$ itself is compact. Such equivalence relations include
free group actions as above by considering the orbit relation
\[
R_{\varphi} = \{ (x, \varphi^{g}(x)) \mid x \in X, g \in G \}.
\]
The map from $X \times G$ to $R_{\varphi}$ sending a pair $(x, g)$ to $(x, \varphi^{g}(x))$
is surjective (by the definition of $R_{\varphi}$) and injective (by freeness). The topology 
on  $R_{\varphi}$ is obtained by simply transferring the product topology ($G$ considered as a discrete space)
from $X \times G$.

This extension of the class of systems we are considering, 
from group actions to \'{e}tale equivalence relations,
  is an important one: it admits
 so-called AF-equivalence relations. Briefly, an \'{e}tale equivalence relation
$R$ is AF (or approximately finite) if it can be written as an increasing sequence of 
compact, open subequivalence relations. Such relations have a presentation by 
means of a combinatorial object - a Bratteli diagram. This class is at once rich, but also tractable.

If $R$ is an equivalence relation on $X$ (even without $R$ having any topology itself), 
we say that such a  system $(X,R)$ is minimal if every $R$-equivalence class is dense in $X$.
In addition, if $R$ is an \'{e}tale  equivalence relation on $X$, we say
that a probability measure $\mu$ on $X$ is $R$-invariant if,  $E$ and  $F$ are Borel subsets of $X$ such that there exists a Borel bijection $f:E \rightarrow F$ contained in $R$, then $\mu(E) = \mu(F)$. We refer the reader to \cite{Ren:book,PPZ:survey,GPS:affable}.
We let $M(X,R)$ denote the set of $R$-invariant probability measures on $X$. This is a
weak* compact, convex set, in fact, a Choquet simplex and is non-empty  whenever $R$ is an AF-relation or arises from a free action of an 
amenable group. We say that $(X,R)$ is \emph{uniquely ergodic} if the set $M(X,R)$ has exactly one element.

We introduce an invariant which is an ordered abelian group. This first appeared in 
\cite{GPS:Zaction} (in the case of $\Z$-actions) as the quotient of a K-theory group by its subgroup of infinitesimal elements. Both of these ordered groups were shown to be dimension groups in the case of $\Z$-actions \cite{HPS:Cantor,GPS:Zaction}.
In fact, we do not need this interpretation of the group here, so we provide a less technical definition (which is equivalent in the earlier case). We let $C(X, \Z)$ denote the set of continuous, integer-valued
functions on $X$. It is an abelian group with the operation of pointwise addition.
We let $B_{m}(X,R)$ denote the subgroup of all functions $f$ such that $\int_{X} f d\mu =0$, 
for all $\mu$ in $M(X,R)$. The quotient group $C(X, \Z)/B_{m}(X,R)$ is denoted $D_{m}(X,R)$.
(Here, the $m$ is used to suggest `measure'.) For a function $f$ in $C(X,\Z)$, we denote 
its class in the quotient by $[f]$. Of course, this is a countable abelian group, but
it is also given an order structure \cite{GPS:Zaction} by defining the positive cone $D_{m}(X,R)^{+}$
as the set of all $[f]$, where $ f \geq 0$. It also has a distinguished positive element,
$[1]$, where $1$ denotes the constant function with value $1$.
Our invariant is the triple $(D_{m}(X,R), D_{m}(X,R)^{+}, [1])$.

Our primary interest is in the notion of orbit equivalence.

\begin{defn}
Let $X$ and $X'$ be two topological spaces and let $R$ and $R'$ be equivalence relations
on $X$, $X'$, respectively. We say that $(X,R)$ and $(X', R')$ are \emph{orbit equivalent }
if there is a homeomorphism $h : X \rightarrow X'$ such that $h \times h(R) = R'$.
\end{defn}

The overall objective (well out of reach at this point) is to try to classify such dynamical systems up to orbit equivalence. This is the topological analogue of an important program in measurable dynamics initiated by Henry Dye \cite{D:orbit} and continued by many others, notably Ornstein and Weiss \cite{OW:shortamen,OW:amenable} and Connes, Feldman and Weiss \cite{CFW:amenable}. There has also been considerable work done in the Borel category, for example see \cite{JKL:Borel}.

The first observation is the following, which is a fairly simple consequence of
the definitions (see \cite{GPS:affable,PPZ:survey}).

\begin{thm}
\label{thm:Dm_invariant}
Let $X$ be a Cantor set and let $R$ be an equivalence relation with an \'{e}tale topology. 
The group, with positive cone and distinguished positive element, 
$(D_{m}(X, R), D_{m}(X, R)^{+}, [1]) $, is an invariant of orbit equivalence.
\end{thm}

As we will see in a moment, the AF-relations play a distinguished r\^{o}le, and so we make the following
definition.

\begin{defn}
An equivalence relation $R$ on $X$ is \emph{affable} if it is orbit equivalent to an AF-relation.
\end{defn}

The term comes from the fact that, if $h$ is an orbit equivalence between two  systems, $(X,R)$ and
$(X',R')$, where $R$ and $R'$ have \'{e}tale topologies, the map $h \times h : R \rightarrow R'$
may not be a homeomorphism. For example, if $\varphi$ is a minimal  action of $\Z$ on $X$, then
it is shown in \cite{GPS:Zaction} that $R_{\varphi}$ is orbit equivalent to an AF-relation. In this case, the map
$h \times h$ cannot be a homeomorphism between the two relations with their \'{e}tale topologies,
since $R_{\varphi}$ is generated, as an equivalence relation, by $\{ (x, \varphi^{1}(x)) \mid x \in X \}$,
which is compact, while it is easy to see from the definition of  an AF-relation that it cannot be generated
by any compact subset. In this case, we can regard the map $h \times h$ as giving a new \'{e}tale topology
on $R_{\varphi}$ in which it is AF. That is, it is AF-able or affable.
The following result is proved in \cite{GPS:Zaction}.

\begin{thm}
\label{thm:1_AF_Z}
Let $(X,R)$ and $(X',R')$ be two minimal equivalence relations on Cantor sets  which are either 
AF-relations or arise from actions of the group $\Z$. Then they are orbit equivalent if and only if
\[
(D_{m}(X, R), D_{m}(X, R)^{+}, [1]) \cong (D_{m}(X', R'), D_{m}(X', R')^{+}, [1]),
\]
meaning that there is a group isomorphism between $D_{m}(X, R)$ and $D_{m}(X', R')$ which is a bijection between
positive cones and preserves the class of $1$.
\end{thm}

In order to extend this result to more general group actions, it suffices to show that, given
a minimal action, $\varphi$, of some group, $G$, on $X$, the orbit relation $R_{\varphi}$ is orbit 
equivalent to an AF-relation; i.e. is affable. The aim of this paper is to establish this for
$G = \Z^{2}$. Our main result is the following.

\begin{thm}
\label{thm:1_affable}
Let $\varphi$ be a free, minimal action of $\Z^{2}$ on the Cantor set $X$. Then the
equivalence relation $(X, R_{\varphi})$ is affable.
\end{thm}

The proof is quite long and it will occupy the rest of the paper. Assuming this for the moment,
it immediately gives the following extension of the earlier result.

\begin{thm}
\label{thm:1_main}
Let $(X,R)$ and $(X',R')$ be two minimal equivalence relations on Cantor sets  which are either 
AF-relations or arise from free actions of the group $\Z$ or of the group $\Z^{2}$. 
Then they are orbit equivalent if and only if
\[
(D_{m}(X, R), D_{m}(X, R)^{+}, [1]) \cong (D_{m}(X', R'), D_{m}(X', R')^{+}, [1]),
\]
meaning that there is a group isomorphism between $D_{m}(X, R)$ and $D_{m}(X', R')$ which is a bijection between
positive cones and preserves the class of $1$.
\end{thm}

Let us add some remarks on the range of the invariant and derive one interesting consequence. The range of the invariant, $D_{m}$, for minimal AF-relations is precisely the collection of simple, acyclic dimension groups with no non-trivial infinitesimal elements. In \cite{HPS:Cantor} and \cite{GPS:Zaction}, it is shown that the range of the invariant is exactly the same for minimal $\Z$-actions on the Cantor set. It follows from  Theorems \ref{thm:Dm_invariant} and \ref{thm:1_main} that the range for minimal $\Z^{2}$-actions is contained in this same collection. At this point, we do not have an exact description of this range. However, it does  follow that every minimal, free $\Z^{2}$-action on a Cantor set is orbit equivalent to a $\Z$-action.

The following two corollaries are immediate consequences of the main theorem and the definitions. We omit the proofs.

\begin{cor}
 \label{cor:1_main_measures}
Let $(X,R)$ and $(X',R')$ be two minimal equivalence relations on Cantor sets  which are either 
AF-relations or arise from free actions of the group $\Z$ or of the group $\Z^{2}$. 
Then they are orbit equivalent if and only if there exists a homeomorphism
$h: X \rightarrow X'$ which implements a bijection between the sets $M(X,R)$ and 
$M(X',R')$.
\end{cor}

\begin{cor}
 \label{cor:1_main_unique_ergodic}
Let $(X,R)$ and $(X',R')$ be two minimal, uniquely ergodic equivalence relations on Cantor sets  which are either 
AF-relations or arise from free actions of the group $\Z$ or of the group $\Z^{2}$. 
Suppose that $M(X,R) = \{ \mu \}$ and 
 $M(X',R') = \{ \mu' \}$. Then the two systems are orbit equivalent if and only if
\[
 \{ \mu(U) \mid U \subset X, U \text{ clopen } \} = 
\{ \mu'(U') \mid U' \subset X', U' \text{ clopen } \}.
\]
\end{cor}

In an earlier, unpublished paper, three of the authors of this paper gave a proof of the main result under the additional hypothesis that the $\Z^{2}$-action had sufficiently many `small, positive' cocycles. See \cite{GPS:abel}. While we still believe that the issue of the existence of such cocycles is an important one, it does not play a r\^{o}le in this paper. The second author, building on some ideas in the unpublished work, gave  proofs of the main result in two special cases in \cite{Ma:Cantor,Ma:tiling}. The current paper is a result of the synthesis of ideas in these two papers and the earlier unpublished one.

We also mention that the result also holds without the hypothesis of freeness of the action, as follows. Let $x$ be in $X$ and let $H_{x} = \{ n \in \Z^{2} \mid \varphi^{n}(x) = x \}$. It is easy to see that if $x$ and $x'$ are in the same orbit, then $H_{x} = H_{x'}$. Secondly, if $x_{i}, i \geq 1$ is a sequence converging to some point $x$ in $X$ and $n \in H_{x_{i}}$, for all $i$, then $n \in H_{x}$ also. From these facts and the minimality of the action, it can be shown that $H_{x}$ is the same for all $x$ in $X$. So the orbits may also be realized as a free action of the group $\Z^{2}/H$. The only possibilities for this quotient group are finite groups (which cannot act minimally on an infinite space), $\Z$, $\Z^{2}$ and $\Z \oplus (\Z / m \Z)$, for some $m \geq 2$. In the last case, the result is due to  Johansen \cite{J:thesis}.

The paper is organized as follows. The second section contains preliminary material. The most important item is the absorption theorem and a result, Theorem \ref{thm:affapp}, which allows us to deduce Theorem  \ref{thm:1_affable}. The rest of the paper is devoted to proving Theorem \ref{thm:affapp}. The following two sections
describe the construction of special tessellations of the plane from minimal actions
of $\Z^{2}$. The fifth section describes a process for refining these tessellations and the final section
uses these tools in giving a proof of Theorem \ref{thm:affapp}.

The authors would like to thank the referees for many helpful suggestions.

\section{Preliminaries}

Throughout this section, $X$ will denote a compact, totally disconnected metrizable space.
The following definitions and results are all taken from \cite{GMPS:affable}, or an earlier version \cite{GPS:affable}. Our main goal is to describe the absorption theorem of \cite{GMPS:affable} and state a theorem from which our main result, Theorem \ref{thm:1_affable} may be deduced.

If $R$ and $S$ are two equivalence relations on a set $X$, we define
\[
R \times_{X} S = \{ ((x,y), (y,z)) \mid (x,y) \in R, (y,z) \in S \}
\]
and we define $r, s: R \times_{X} S \rightarrow X$ by
\[
r((x,y),(y,z)) = x, s((x,y),(y,z)) = z.
\]
 If $R$ and $S$  have topologies, we give
$R \times_{X} S $ the relative topology from the product.

If $R$ is an equivalence relation on $X$ and $Y$ is any subset of $X$, we let
$R|Y$ denote the restriction of $R$ to $Y$, that is, $R|Y = R \cap (Y \times Y)$.
Moreover, we say that $Y$ is \'{e}tale in $(X, R)$ if it is closed and $R|Y$, with its relative topology from $R$, is \'{e}tale \cite{GPS:affable}.

\begin{defn}
Let $(X,R)$ be an \'{e}tale equivalence relation.
 A compact \'{e}tale equivalence relation $K$ on
$X$ is \emph{transverse} to $R$ if
\begin{enumerate}
\item
$K \cap R = \Delta_{X} = \{ (x,x) \mid x \in X \}$,
\item
There is a homeomorphism, $h: R \times_{X} K \rightarrow K  \times_{X} R$ such that
$r \circ h = r$ and $s \circ h = s$.
\end{enumerate}
\end{defn}

One simple, but important class of examples is the following. Suppose that $(X,R)$
is an \'{e}tale equivalence relation and $\alpha$ is an action of the finite group
$G$ on $X$ such that 
\begin{enumerate}
\item
$\alpha_{g} \times \alpha_{g}(R) = R$, for all $g$ in $G$,
\item
$\alpha_{g} \times \alpha_{g}: R \rightarrow R$ is a homeomorphism, for all $g$ in $G$,
and
\item
$(x, \alpha_{g}(x))$ is not in $R$, for any $x$ in $X$, $g \neq e$ in $G$ (in particular, the action is free).
\end{enumerate}
Then the relation $K = \{(x, \alpha_{g}(x)) \mid x \in X, g \in G \}$ is transverse
to $R$ via the map $h((x,y),(y, \alpha_{g}(y))) = ((x, \alpha_{g}(x)), (\alpha_{g}(x), \alpha_{g}(y)))$, for $(x,y)$ in $R$ and $g$ in $G$.

We state the following two results from \cite{GMPS:affable} for completeness.

\begin{thm}
\label{thm:2_transverse}
Let $R$ be an \'{e}tale equivalence relation and $K$ be a compact \'{e}tale 
equivalence relation on
$X$ which is transverse to $R$. The equivalence relation generated by $R$ 
and $K$, denoted $R \vee K$ 
is equal to $r \times s (R \times_{X} K)$. More precisely, the map 
from $R \times_{X} K$ to $R \vee K$ defined by
sending $((x,y),(y,z))$ to $(x,z)$ is a bijection, and with the topology 
from $R \times_{X} K$ transferred
by this map, $R \vee K$ is \'{e}tale. Moreover, this is the unique 
\'{e}tale topology which extends that
of $R$ and of $K$.
\end{thm}

\begin{thm}
\label{thm:2_transverse_AF}
Let $R$ be an AF-relation and $K$ be a compact \'{e}tale 
equivalence relation on
$X$ which is transverse to $R$. Then with the unique topology 
extending that of $R$ and $K$, $R \vee K$
is an AF-relation.
\end{thm}

With this terminology and notation, the absorption theorem is the following. It will be  the essential technical tool in the
proof of the main result. We refer the reader to \cite{GMPS:affable} for its proof and further discussion.

\begin{thm}
\label{thm:2_absorb}
Let $(X,R)$ be a minimal AF-relation, let $Y$ be a closed subset of $X$ and let $K$
be a compact \'{e}tale equivalence relation on $Y$. Suppose that the following hold:
\begin{enumerate}
\item
$\mu(Y) = 0$, for all $R$-invariant probability measures $\mu$ on $X$,
\item
$Y$ is an \'{e}tale subset of $(X,R)$,
\item
$K$ is transverse to $R | Y$.
\end{enumerate}
Then there is a homeomorphism $h$ of $X$ such that
\begin{enumerate}
\item
\[
h \times h (R \vee K) = R,
\]
where $R \vee K$ is the equivalence relation generated by $R$ and $K$, and 
\item
$h(Y)$ is \'{e}tale in $(X,R)$ and $\mu(h(Y)) = 0$, for all $R$-invariant probability measures $\mu$ on $X$,
\item
\[
h|_{Y} \times h |_{Y}: (R | Y)  \vee K \rightarrow R | h(Y)
\]
is a homeomorphism.
\end{enumerate}
In particular, $R \vee K$ is affable.
\end{thm}

We will apply this result to our minimal free $\Z^{2}$-actions in the following fashion.

\begin{thm}
 \label{thm:affapp}
Let $\varphi$ be a free minimal action of $\Z^{2}$ on a Cantor set $X$. Then there exist:
\begin{enumerate}
 \item a minimal AF-relation $R$ which is an open subrelation of $R_{\varphi}$,
\item a closed set $B^{2} \subset X$ which is  \'{e}tale in $(X, R)$ and such that $\mu(B^{2}) = 0$, for all $R$-invariant probability measures $\mu$,
\item a compact \'{e}tale equivalence relation $K^{2}$ on $B^{2}$ which is transverse to $R | B^{2}$,
\item a closed set $B^{3} \subset B^{2}$ which is \'{e}tale in $(B^{2}, (R | B^{2}) \vee K^{2})$,
\item a compact \'{e}tale equivalence relation $K^{3}$ on $B^{3}$ which is transverse to $((R | B^{2}) \vee K^{2}) | B^{3}$.
\end{enumerate}
Moreover, we have $R_{\varphi} = R \vee K^{2} \vee K^{3}$.
\end{thm}

The proof of this result will occupy the rest of the paper. We now prove that our main result, Theorem \ref{thm:1_affable}, follows from this.

 Given a free minimal action $\varphi$ 
of $\Z^{2}$ on $X$, we construct $R, B^{2}, K^{2}, B^{3}, K^{3}$ of Theorem 
\ref{thm:affapp}. 
We apply the absorption Theorem \ref{thm:2_absorb}
using the  AF-relation $R$ on $X$, the closed set $Y = B^{2}$ and the equivalence relation
$K = K^{2}$ on $B^{2}$. The hypotheses are satisfied by 
\ref{prop:6_B2_etale}, \ref{prop:6_K2_transverse}, \ref{prop:6_R_minimal} and \ref{prop:6_B2_measure}. 
We obtain a homeomorphism $h_{2}$ of $X$ such that $h_{2} \times h_{2}(R \vee K^{2}) = R$,
$h_{2}(B^{2})$ is \'{e}tale for $R$ and has measure zero, for all $R$-invariant measures and
\[
h_{2}|_{B^{2}} \times h_{2}|_{B^{2}}: (R|B^{2}) \vee K^{2} \rightarrow R|h_{2}(B^{2})
\]
is a homeomorphism. 

Next, we apply the absorption theorem a second time with
the same $R$, $Y = h_{2}(B^{3}), K= h_{2} \times h_{2}(K^{3})$. That the appropriate hypotheses
hold is an immediate consequence of \ref{prop:6_B3_etale}, \ref{prop:6_K3_transverse}  
and the conclusions from the the first application of the absorption theorem. We obtain
a homeomorphism $h_{3}$ of $X$ such that 
\[
h_{3} \times h_{3}(R \vee (h_{2} \times h_{2})(K^{3})) = R.
\]
It follows that $h_{3} \circ h_{2}$ is a homeomorphism of $X$ and, by Proposition \ref{prop:6_generate},
\begin{eqnarray*}
(h_{3} \circ h_{2}) \times (h_{3} \circ h_{2}) (R_{\varphi}) & = & 
                 (h_{3} \circ h_{2}) \times (h_{3} \circ h_{2}) (R \vee K^{2} \vee K^{3}) \\
     &  =   &  (h_{3} \times h_{3}) \circ (h_{2} \times h_{2}) (R \vee K^{2} \vee K^{3}) \\
    &  =  &  (h_{3} \times h_{3})( (h_{2} \times h_{2})(R \vee K^{2}) \vee ((h_{2} \times h_{2})(K^{3})) \\
   &  =  &    (h_{3} \times h_{3})(R \vee ((h_{2} \times h_{2})(K^{3})) \\
   &  =  &  R.
\end{eqnarray*}
This completes the proof of Theorem \ref{thm:1_affable}.

\section{Tessellations of $\R^{2}$}

In this section, we make preparations for our proof of Theorem \ref{thm:affapp}.
We will be constructing various tessellations of $\R^{2}$ (the definition
will follow below).  The first tool in this is the notion of Voronoi tessellation. For our purposes, Voronoi tessellations have three drawbacks. We will want the vertices of our tessellations to be the intersections of precisely three tiles; generically, this is the case for Voronoi cells, but it can occur that more than three cells meet at a point. This problem is relatively easy to handle; we consider the dual tessellation and sub-divide until it is a triangulation. The second problem concerns the combinatorics of our triangulation. Briefly, we require that if there are three cells, any two of which meet, then all three should have non-empty intersection. This problem is solved by removing triangles which are contained in others. The third problem is more serious. We would like to know that cells which are disjoint should be separated in some controlled manner. We are able to do this by moving the vertices of the Voronoi tessellation to the incentres of the triangles of the
  dual. (There are undoubtedly other ways of doing this, but this seems the most convenient.)  

These second and third difficulties are the points in our arguments which do \emph{not} seem to generalize to the situation of $\R^{d}$, for $d > 2$. Ultimately, we introduce the notion of a \emph{simplicial} tessellation. Here, the combinatorics of the cells is carefully controlled. This allows to give combinatorial rather than geometric descriptions of edges and vertices.

Other constructions of this type have been made in the case of Cantor minimal systems; first by Forrest \cite{For:Z2action}, also described in \cite{Ph:Z2action} and more recently by Lightwood and Ormes \cite{LO:Z2action}.

Let $d(\cdot, \cdot)$ denote the usual metric on $\R^{2}$.
For any non-empty set $A$ and $u$ in $\R^{2}$, we let $d(u,A) = \inf \{ d(u,v) \mid v \in A \}$
and for two non-empty  sets, $A_{1}, A_{2}$, we let 
$d(A_{1}, A_{2}) = \inf \{ d(u,v) \mid u \in A_{1}, v \in A_{2} \}$. We let $B(u,r)$ denote 
the open ball of radius $r > 0$ centred at $u$ in $\R^{2}$ and $S(u,r)$ denote the sphere
of radius $r$ centred at $u$.
Given two distinct points $u$ and $v$ in $\R^{2}$, we let
\[
\overline{uv} = \{ tu + (1-t)v \mid 0 \leq t \leq 1 \}
\]
be the closed line segment connecting $u$ and $v$.

We begin by establishing some basic geometric facts in the plane.

\begin{lemma}
\label{lemma:3_angle}
Let $u_{1}, u_{2}, u_{3}$ be three distinct points on $S(0, r)$ in $\R^{2}$. If $r \leq 2M$ and 
$  d(u_{i}, u_{j}) \geq M$, for all $i \neq j$, then there are constants, $0 < \alpha, \beta < \pi$,
independent of $M$, such that 
\[
\alpha \leq \angle u_{1}u_{2}u_{3} \leq \beta.
\]
\end{lemma}

\begin{proof}By the law of sines, we have 
\[
 \frac{d(u_{1},u_{3})}{ \sin (\angle u_{1}u_{2}u_{3}) } = 2r.
\]
It follows that
\[
 \sin (\angle u_{1}u_{2}u_{3})  = \frac{d(u_{1}, u_{3})}{2r} \geq \frac{M}{4M} = \frac{1}{4},
\]
and so $\angle u_{1}u_{2}u_{3} \geq \arcsin(1/4)$. On the other hand, we also have
\begin{eqnarray*}
\angle u_{1}u_{2}u_{3} & = & \pi - \angle u_{2}u_{3}u_{1} - \angle u_{3}u_{1}u_{2} \\
                 & \leq  & 2(\frac{\pi}{2} - \arcsin(1/4)) \\
      & =  & 2 \arccos( 1/4).
\end{eqnarray*}
\end{proof}

Recall that the three angle bisectors of a triangle meet at a point called the incentre of the triangle.

\begin{lemma}
\label{lemma:3_incentre}
Given $0 < \alpha, \beta < \pi$, there exists a positive constant $b$ satisfying the following.
Let $t = \bigtriangleup u v w$ be a triangle such that  all side lengths are at least $M$
and all angles are between $\alpha$ and $\beta$. 
Let $c(t)$ be the incentre of the triangle $t$. Then 
\begin{enumerate}
\item
 $B(c(t), bM) \subset t$. 
\item
for  two triangles, $t, t'$,  as above which share an edge, but have disjoint interiors,
the line segment $\overline{c(t)c(t')}$ is contained in $t \cup t'$ and meets the common edge.
\item
for  $t, t'$ as above and supposing that $\overline{uv}$ is the common edge, we have 
\[
d( \{u, v \}, \overline{c(t)c(t')}) \geq bM.
\]
\end{enumerate}
\end{lemma}

\begin{proof}

We let 
\[
b = \frac{ \cos^{2}( \frac{\beta}{2}) \sin( \frac{\alpha}{2})}{2}.
\]

For the first part, let us assume for convenience that $v$ is the origin and the angle bisector of
$\angle u v w$ is the $x$-axis, with $u$ in the first quadrant and $w$ in the fourth quadrant.
We may also assume that the $x$-coordinate of $u$ is less than or equal to that of $w$. 
The complement of the triangle (in the right half-plane) is covered by three sets: the part of 
the half-plane above the line through $v$ and $u$, the part below the line through $v$ and $w$
and the half-plane to the right of the vertical line passing through $u$. From the lower bound on the
angle $\angle u v w$, the first part is contained
in the part of the right half plane above the line through the origin with angle $\frac{\alpha}{2}$ to
the $x$-axis. Similarly, the third part is contained below the line making angle
$-\frac{\alpha}{2}$ with the axis. The minimum value of the $x$-coordinate of $u$ occurs
when the side length $\overline{uv}$ is as small as possible, namely, $M$, and when the angle of $\overline{uv}$ with
the $x$-axis is as large as possible, namely, $\frac{\beta}{2}$. The minimum value of this
$x$-coordinate is $M \cos(\frac{\beta}{2}) $. From this we see that the triangle contains the
isosceles triangle with vertex at the origin, one side along the vertical line with $x$-coordinate
$M \cos(\frac{\beta}{2})$ and other sides along the lines making angles $\pm \frac{\alpha}{2}$ with the $x$-axis.
This triangle contains 
$B((\frac{M\cos(\frac{\beta}{2})}{2}, 0),\frac{ M\cos(\frac{\beta}{2}) \sin{\alpha/2}}{2})$, 
which in turn contains a ball of radius $bM$, since $b \leq \frac{\cos(\frac{\beta}{2}) \sin(\alpha/2)}{2} $ . 
The incentre of a triangle is the
centre of  the largest ball contained in the triangle and so 
\[
b = \frac{\cos^{2}( \frac{\beta}{2}) \sin( \frac{\alpha}{2})}{2}
\]
satisfies the first condition.

For the second  and third part, assume that the triangle $ t = \bigtriangleup u v w$ has $u$ 
at the origin and $v$ on the positive $x$-axis, 
with $\overline{uv}$ the common edge between the two triangles. Consider the minimum possible value for
the $x$ coordinate of $c(t)$. Recalling that the incentre lies on the angle bisectors
of the triangle, the maximum value of $\angle c(t) u v$ is $\frac{\beta}{2}$ and the minimum
distance from $u$ to $c(t)$, as shown above, is $ \frac{ M\cos(\frac{\beta}{2}) \sin{\alpha/2}}{2} $. 
Hence the minimum value for the $x$-coordinate of $c(t)$ is
$ \cos(\beta/2) \frac{ M\cos(\frac{\beta}{2}) \sin{\alpha/2}}{2} = bM$. 
The second and third parts follow at once.
\end{proof}

A \emph{tessellation} of $\R^{2}$,
 $\mathcal{T}$, 
 is a collection of closed polygons which cover $\R^{2}$, with pairwise disjoint
interiors. We also let 
$\mathcal{T}^{1}$ denote the   edges which form their boundaries  
and $\mathcal{T}^{0}$ denote the
points which are the vertices.  For the moment, the edges will be line segments between vertices, but later, we will allow an edge to be the union of a finite collection of line segments, homeomorphic to the unit interval. Here the vertices will mean the endpoints of the edges, and not the endpoints of the individual line segments making up the edges.

Let $P$ be a countable subset of $\R^{2}$. For a real number $M \geq 0$, we say that $P$ is 
$M$-separated if $d(u,v) \geq M$, for all $ u \neq v$ in $P$. We also say that $P$ is
$M$-syndetic if $ \cup_{u \in P} B(u, M) = \R^{2}$.

We will show how to start from an $M$-separated, $2M$-syndetic set $P$ and construct
a tessellation $\mathcal{T}_{P}$ with various nice properties. Before beginning, 
we make the following remark. In the actual application, we will begin
with a $\varphi$-regular collection (the definition will be given in the following section, but it is not needed for the moment),
 $P(x), x\in X$, and construct a collection
of tessellations, $\mathcal{T}_{P(x)}$. It is worth noting as we proceed, that
all of our constructions are `locally derived' in the  sense that there is an absolute constant $R$ so that the vertices, edges and so on in some region $A$, depend only on the data in the ball around $A$ of radius $R$. By application of \ref{prop:4_local_derive}, the collection, $\mathcal{T}_{P(x)}, x \in X$, 
is $\varphi$-regular.

Fix $M > 0$ and let $P$ be an $M$-separated, $2M$-syndetic  subset of $\R^{2}$. 
We first define the
associated Voronoi tessellation of $\R^{2}$, $\mathcal{T}_{v}$ as follows. 

 For each $u$ in $P$, let
\[
T(u) = \{ v \in \R^{2} \mid d(v,u) \leq d(v,P) \},
\]
which is a polygon with $u$ in its interior.
$\mathcal{T}_{v}(P)$ denotes the collection of $T(u), u \in P$.
For any two elements $u, v$ of $P$, we define the edge separating $u$ and $v$ by
\[
u|v = \{ w \in \R^{2} \mid d(w,u) = d(w,v) = d(w, P) \}
\]
provided that set is infinite.
The edges of the Voronoi tessellation are the infinite sets 
$u|v$, $u, v \in P$ and $\mathcal{T}_{v}^{1}(P)$
denotes the edges. The
vertices are the set of all $w$ such that $S(w, d(w,P))$ contains at least three 
points of $P$ and $\mathcal{T}_{v}^{0}(P)$ denotes the vertices.
 For such a vertex $w$, we let $P(w)$ denote $ S(w, d(w,P)) \cap P$.
We adopt the convention that we list these points as $\{ u_{1}, u_{2}, \ldots, u_{k} \}$
 in counter-clockwise order around $w$, beginning at the right horizontal. That is,  for all
$1 \leq i \leq k$, $u_{i} -w = (r \cos(\theta_{i}), r \sin( \theta_{i}) )$, with
$0 \leq \theta_{1} < \theta_{2} < \cdots < \theta_{k} < 2\pi$. 

We describe the dual Voronoi tessellation of $P$. The vertices, $\mathcal{T}_{dv}^{0}(P)$,
 are the points of $P$. For each
pair of $u, v$ in $P$ such that $u|v$ is in $\mathcal{T}_{v}^{1}(P)$, the line segment $\overline{uv}$
is an edge and we let $\mathcal{T}_{dv}^{1}(P)$ denote the set of edges in the dual. Finally, the
faces of the dual, which we denote by $\mathcal{T}_{dv}^{2}(P)$ are indexed by the vertices in 
$\mathcal{T}_{v}^{0}(P)$: for each such $w$, the associated face of the dual is the polygon 
with vertices $P(w)$.

Generically, the faces of the dual Voronoi tessellation are triangles. However, it is possible that
some faces have more than three edges. We want to correct this, producing a triangulation
of the plane which we denote by $\mathcal{T}_{dt}(P)$, called 
the dual Voronoi triangulation associated with $P$. It has the same
vertex set, namely $P$. The edge set, $\mathcal{T}_{dt}^{1}(P)$ contains $\mathcal{T}_{dv}^{1}(P)$ and, 
for every $w$ in $\mathcal{T}_{v}^{0}(P)$ with $P(w) = \{ u_{1}, u_{2}, \ldots, u_{k} \}$ with $k \geq 3$,
we also include the edges $\overline{u_{1}u_{i}}$, with $2 < i < k$. The faces, denoted
$\mathcal{T}_{dt}(P)$, are then all the triangles
$\bigtriangleup u_{1}u_{i}u_{i+1}$, $1 < i < k$. Notice that in the case that $k=3$, this just yields the same
(single) face as before.

For each triangle in $\mathcal{T}_{dt}(P)$, the side lengths are at least $M$ and its circumcircle has radius less than $2M$ because $P$ is $M$-separated and $2M$-syndetic. It follows from Lemma \ref{lemma:3_angle} that all angles of any triangle in  $\mathcal{T}_{dt}(P)$ are between $\alpha$ and $\beta$.

Next, we modify $\mathcal{T}_{dt}(P)$ further and define a new triangulation 
$\mathcal{T}_{*}(P)$ as follows. Let $F$ denote the set of all vertices of $\mathcal{T}_{dt}(P)$ which
are contained in the interior of some triangle whose edges are in $\mathcal{T}_{dt}^{1}(P)$ (although the triangle itself may be the union of ones in 
$\mathcal{T}_{dt}^{2}(P)$). That is, $F$ is all $u$ in $P$ such that there exist $u_{1}, u_{2}, u_{3}$ in $P$ with $\overline{u_{1}u_{2}}$,
$\overline{u_{1}u_{3}}$, $\overline{u_{2}u_{3}}$ in $\mathcal{T}_{dt}^{1}(P)$ and 
$u$ in the interior of $\triangle u_{1}u_{2}u_{3}$. We remove such vertices from 
$\mathcal{T}_{dt}^{0}(P)$ and define 
\begin{eqnarray*}
 \mathcal{T}_{*}^{0}(P) & = & P \setminus F,  \\
\mathcal{T}_{*}^{1}(P) & = & \{ \overline{uv} \in \mathcal{T}_{dt}^{1}(P) \mid 
                                            u, v \in  \mathcal{T}_{*}^{0}(P) \}, \\ 
\mathcal{T}_{*}^{2}(P) & = & \{ \triangle u_{1}u_{2}u_{3} \mid 
                   \overline{u_{1}u_{2}}, \overline{u_{1}u_{3}}, \overline{u_{2}u_{3}}
                       \in \mathcal{T}_{*}^{1}(P) \}.
\end{eqnarray*}
In simpler terms, we consider all triples of vertices $u,v,w$ such that each pair is joined by a line segment. Such triples may be ordered: one is smaller if the triangle which it determines is contained in the other. Our triangulation $\mathcal{T}_{dt}$ is given by the minimal triangles. By the separation properties of $P$, it is clear that any chain in our ordered set of triples must be finite. We consider only the maximal triangles and this forms $\mathcal{T}_{*}(P)$. It is clear that $\mathcal{T}_{*}(P)$ is a triangulation of  $\R^{2}$. Moreover, if $u,v,w$ are vertices such that any two are the endpoints of an edge, then the interior of the triangle they determine is an element of $\mathcal{T}_{*}(P)$.

Any triangle in $\mathcal{T}_{*}(P)$ is a union of triangles in $\mathcal{T}_{dt}(P)$ and it follows that the  angle of any triangle is bounded below by $\alpha$. From this, we can conclude that any such angle is bounded above by $\beta = \pi - 2\alpha$.  Therefore, Lemma \ref{lemma:3_incentre} may be applied to triangles in $\mathcal{T}_{*}(P)$.

If $u, v$ are in $P$ and $\overline{uv}$ is an edge of $\mathcal{T}_{dt}(P)$, then the Voronoi cells at $u$ and $v$ meet.
Since $P$ is $2M$-syndetic, the point where these cells meet is distance less than $2M$ from $u$ and $v$. We conclude that the length of the edge $\overline{uv}$ is less than $4M$. It then follows that any edge of $\mathcal{T}_{*}(P)$ is also of length less than $4M$. 

Using this triangulation, $\mathcal{T}_{*}(P)$,  we define a new tessellation, denoted
$\mathcal{T}_{P}$, as follows. The vertex set, $\mathcal{T}_{P}^{0}$, is the collection of incentres, $c(t)$, where $t$ is in $\mathcal{T}_{*}^{2}(P)$. The edge set,
$\mathcal{T}_{P}^{1}$, is all $\overline{c(t)c(t')}$, where $t,t'$ share an edge. The elements of $\mathcal{T}_{P}$ are the polygons with these edges. Each such polygon, $t$, contains a unique point of $P \setminus F$, which we denote by $P(t)$.

One of the nice properties of this tessellation is summarized in the following definition.
However, we should make some remarks, since this definition anticipates some
constructions which occur later. As we modify our tessellations in section 5, we will have geometric objects which cover the plane, but are not precisely polygons. The first problem is that, although the objects may even be polygons, we will prefer to think of an edge as the intersection of two elements of the tessellation, which may, in fact, be a union of line segments. In other terms, we want to have a combinatorial notion of edge and vertex, rather than a geometric one. Another aspect which is allowed by this definition, is that the regions may actually be disconnected unions of polygons.

\begin{defn}
A \emph{simplicial tessellation} is a collection, $\mathcal{T}$, of compact
subsets of $\R^{2}$
 which cover $\R^{2}$, have pairwise disjoint interiors and satisfy
\begin{enumerate}
\item
if $t_{1}, t_{2}, \ldots , t_{k}$ and $t'_{1}, t'_{2}, \ldots , t'_{l}$ are elements such that
\[
\cap_{i=1}^{k} t_{i} = \cap_{j=1}^{l} t'_{j} \neq \emptyset,
\]
then
\[
\{t_{1}, t_{2}, \ldots, t_{k} \} = \{ t'_{1}, t'_{2}, \ldots t'_{l} \},
\]
\item
if $t_{1}, t_{2}, \ldots, t_{k}$ are elements such that for all $i,j$, 
$t_{i} \cap t_{j} \neq \emptyset$, 
then 
\[
\cap_{i=1}^{k} t_{i} \neq \emptyset.
\]
\end{enumerate}
\end{defn}

Note that in the following example, $t_{1}, t_{2}, t_{3}$
 fail to satisfy the last condition.

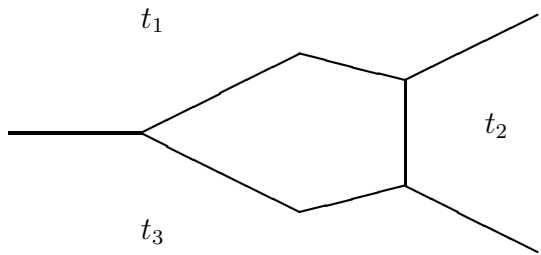
\begin{picture}(200,180)(0,0)

\thicklines

\put(0,90){\line(1,0){50}}
\put(50,90){\line(2,1){60}}
\put(50,90){\line(2,-1){60}}
\put(110,120){\line(4,-1){40}}
\put(110,60){\line(4,1){40}}
\put(150,70){\line(0,1){40}}
\put(150,70){\line(2,-1){50}}
\put(150,110){\line(2,1){50}}

\thinlines

\put(50,130){$t_{1}$}
\put(180,90){$t_{2}$}
\put(50,50){$t_{3}$}

\end{picture}

\begin{defn}
Let $\mathcal{T}$ be any tessellation of $\R^{2}$.
We say that $\mathcal{T}$ has \emph{capacity} $C > 0$ if each element of $\mathcal{T}$
contains a ball of radius $C$. We say that $\mathcal{T}$ is \emph{$K$-separated}, for  $K > 0$ if, for 
any two disjoint elements of $\mathcal{T}$, $t,t'$, we have $d(t, t') \geq K$. Finally, 
we define the \emph{diameter} of $\mathcal{T}$ to be the supremum of the diameter of its elements.
\end{defn}

\begin{prop} 
\label{prop:3_tessellate}
There are  constants, $ b > 0 , E \geq 3, a > 0$, such that, for any
set $P$ contained in $\R^{2}$ which is $M$-separated and $2M$-syndetic, we have the following.
\begin{enumerate}
\item
$\mathcal{T}_{P}$ is a simplicial tessellation,
\item
each element of $\mathcal{T}_{P}$ is a polygon with at most $E$ edges,
\item
$B(P(t), bM)$ is contained in $t$, for any $t$ in $\mathcal{T}_{P}$; in particular, $\mathcal{T}_{P}$ has capacity $bM$, ,
\item
$\mathcal{T}_{P}$ is $bM$-separated,
\item
the angle formed by any two edges of $\mathcal{T}_{P}$ which meet is at least $a$,
\item
for any $u$ in $\mathcal{T}^{0}_{P}$, there exist $t_{1}, t_{2}, t_{3}$ in 
$\mathcal{T}_{P}$ such that $B(u, bM)$ is contained in $t_{1} \cup t_{2} \cup t_{3}$
and, for each $i = 1, 2, 3$, $B(u, bM) \cap t_{i}$ is a sector.

\end{enumerate}
\end{prop}

\begin{proof}
For the first part, begin by considering the intersection of two distinct polygons, $t_{1}, t_{2}$ of $\mathcal{T}_{P}$. It must contain a vertex. Since each vertex is the intersection of three regions, $t_{1}$ and $t_{2}$ must be two of them and hence they share an edge. By Lemma \ref{lemma:3_incentre}, this edge intersects a unique edge of $\mathcal{T}_{*}$, which must be the edge between $P(t_{1})$ and $P(t_{2})$. The first property of a simplicial tessellation follows from this. As for the second property, suppose that $t_{1}, t_{2}$ and $t_{3}$ are three elements, each pair of which meets. We deduce from above that there are edges in $\mathcal{T}_{*}$ between each pair $P(t_{1}), P(t_{2})$ and $P(t_{3})$. But then, it follows from the construction of $\mathcal{T}_{*}$ that the triangle these points determine is in $\mathcal{T}_{*}$. This means that $t_{1}, t_{2}, t_{3}$ have non-empty intersection.

For the second part, we first consider a point $u$ in $P$. Any edge of its Voronoi cell is shared with that another point in $P$, which must be within distance less than $4M$ from $u$. The set $B(u,4M)$ may be covered by a finite number of balls of radius $M/2$, each of which contains at most one point of $P$. So we conclude that the Voronoi cell for $u$ has at most $E_{0}$ edges, for some $E_{0}$. Next, suppose $v$ is any vertex of a Voronoi cell. The set of $u$ in $P$ such that $v$ is in $T(u)$ is contained in $B(v, 2M)$. Again, this ball can be covered by a finite number of balls of radius $M/2$ and so we see that the number of $u$ in $P$ with $v$ in $T(u)$ is bounded by some constant $D_{0}$. This implies that there are at most $D_{0}$ edges containing $v$. Passing now to the dual Voronoi tessellation, $\mathcal{T}_{dv}(P)$, each vertex is contained in at most $E_{0}$ edges and each polygon has at most $D_{0}$ sides. We next add edges in passing to the dual Voronoi triang
 ulation. The number of added edges at any one vertex is at most $D_{0}-2$ in each polygon containing that vertex. This means that a vertex contains at most $E = E_{0} + (D_{0}-2) E_{0} = E_{0}(D_{0}-1)$ edges. The passage to $\mathcal{T}_{*}$ only removes edges, so this bound also holds for $\mathcal{T}_{*}$. This implies that the number of edges around any polygon of $\mathcal{T}_{P}$ is bounded by $E$.

The third statement follows immediately from part 3 of \ref{lemma:3_incentre}. For part 4, consider two
elements $t$ and $t'$ of $\mathcal{T}_{P}$ which do not intersect. First notice
that the minimum distance between two polygons in the plane is achieved when at least one of the points
is at a vertex of its polygon and the other is on an edge. Therefore, it suffices to show the
distance from any vertex of $t$ to $t'$ is at least $bM$. But such a vertex is an incentre for
some triangle in $\mathcal{T}_{*}(P)$ and the conclusion follows at once from part 1
of \ref{lemma:3_incentre}.

For the fifth part, let $t$ be a polygon in $\mathcal{T}_{P}$ and 
consider a vertex of $t$, which is the incentre,  denoted $c(\triangle)$, of some triangle, $\triangle = \triangle uvw$, in 
$\mathcal{T}_{*}(P)$.
$P(t)$ is one of the vertices of this triangle, say $u$. Let $\triangle' = \triangle uvw'$ be another triangle in $\mathcal{T}_{*}(P)$, which shares the edge $\overline{uv}$. Its incentre is another vertex of $t$, connected to $c(\triangle)$ by an edge. Let $\theta = \angle uc(\triangle)c(\triangle')$. We shall show that $\sin{ \theta} \geq \frac{b^{2}}{16}$. Since this holds for any pair of vertices of $t$, we conclude that the angle at any vertex of $t$ is at least $2 \arcsin(\frac{b^{2}}{16})$. Let $x$ be the point where the edge 
$\overline{uv}$ meets $\overline{c(\triangle)c(\triangle')}$, so that
$\theta = \angle uc(\triangle)c(\triangle') = \angle uc(\triangle)x$. 
Consider the area of the triangle $\triangle u c(\triangle) x$. It is clearly equal to 
\[
 \frac{1}{2} d(u, c(\triangle)) d(c(\triangle), x) \sin(\theta).
\]
On the other hand, since the distance from $c( \triangle)$ to the line passing through $u$ and $x$ is at least $bM$ by Lemma \ref{lemma:3_incentre}, this area is greater than or equal to 
\[
 \frac{1}{2} d(u,x) bM.
\]
We have $ d(u, c(\triangle))$ and $ d(c(\triangle), x) $ are at most $4M$ because all the sides of $\triangle$ are less than $4M$. By condition 3 of Lemma \ref{lemma:3_incentre}, $d(u,x)$ is greater than or equal to $bM$. Putting this together, we have
\[
 \sin(\theta) \geq \frac{d(u,x) bM}{d(u, c(\triangle)) d(c(\triangle), x)} \geq 
 \frac{bM \cdot bM}{4M \cdot 4M} =  \frac{b^{2}}{16}
\]
as desired.

The last part is an immediate consequence of part 1 of Lemma \ref{lemma:3_incentre}. 
\end{proof}

\section{$\varphi$-regular tessellations of $\R^{2}$}

The aim of this section is to construct tessellations of $\R^{2}$, satisfying the conditions from the previous section, from a free minimal action of $\Z^{2}$ on a Cantor set.

Let $\varphi$ be a free minimal action of $\Z^{2}$ on a Cantor set $X$.
Suppose that for each $x$ in $X$, we have a subset $P(x)$ of $\R^{2}$. We say this collection
is $\varphi$-regular if
\begin{enumerate}
\item
for any $x$ in $X$ and $n$ in $\Z^{2}$, 
\[
P(\varphi^{n}(x)) = P(x)+ n,
\]
\item
if $x$ is in $X$ and $K \subset \R^{2}$ is compact, 
then there is a neighbourhood  $U$ of $x$ such that 
\[
P(x') \cap K = P(x) \cap K,
\]
for all $x'$ in $U$.
\end{enumerate}

The following result is an easy consequence of the definition and we omit the proof.

\begin{prop}
\label{prop:4_regular}
Let $Y$ be a clopen subset of $X$. The family of sets 
\[
P(x) = \{ n \in \Z^{2} \mid x \in \varphi^{n}(Y) \},
\]
for $x$ in $X$, is $\varphi$-regular.
Conversely, if $P(x)$ is a $\varphi$-regular family, then 
\[
Y = \{ x \in X \mid 0 \in P(x) \}
\]
is clopen.
\end{prop}

We consider a family of tessellations of $\R^{2}$ which are indexed by
the points of $X$, $\mathcal{T}(x)$, $x \in X$. We say that this collection is
$\varphi$-regular if
\begin{enumerate}
\item
for any $x$ in $X$ and $n$ in $\Z^{2}$, 
\[
\mathcal{T}(\varphi^{n}(x)) = \mathcal{T}(x) + n,
\]
\item
if $x$ is in $X$ and $a$ is in $\mathcal{T}(x)$, 
then there is a neighbourhood  $U$ of $x$ such that 
$a$ is in $ \mathcal{T}(x')$, for all $x'$ in $U$.
\end{enumerate}

If $P$ is a $\varphi$-regular family, we 
say it is $M$-syndetic, for some $M \geq 1$,
 ($M$-separated, respectively) if, for each $x$ in $X$, 
$P(x)$ is $M$-syndetic ($M$-separated, respectively).

Let $P(x), P'(x), x \in X$ be two families of subsets of $\R^{2}$. We say that
$P'$ is locally derived from $P$ if there is a constant $R >0$ such that, for any $x_{1}, x_{2}$
in $X$,
and $u_{1}, u_{2}$ in $\R^{2}$, if $u_{1}$ is in $P'(x_{1})$ and 
\[
(P(x_{1}) - u_{1}) \cap B(0,R) = (P(x_{2}) -u_{2}) \cap B(0,R),
\]
then $u_{2}$ is in $P'(x_{2})$. In a similar way, we extend this definition replacing 
either $P$, $P'$ or both with families of tessellations. The following two results are easily derived from the definitions; we omit the proofs.

\begin{prop}
 \label{prop:TP_locally_derived}
Let $P(x), x \in X$ be a $\varphi$-regular family of sets. Then the family
$\mathcal{T}_{P(x)}, x \in X$ of \ref{prop:3_tessellate} is locally derived from $P$.
\end{prop}

\begin{prop}
\label{prop:4_local_derive}
If $P$ is a $\varphi$-regular family and $P'$ is locally derived from $P$, then
$P'$ is also $\varphi$-regular. Analogous statements hold replacing $P$, $P'$ or both
with families of tessellations.
\end{prop}

Finally in this section, we turn to the issue of the existence of $\varphi$-regular, separated and syndetic sets for minimal Cantor $\Z^{2}$-actions. The proof is not new (see \cite{LO:Z2action}, for example), but we provide it here for completeness.

\begin{prop}
\label{prop:4_sep_synd}
Let $(X, \varphi)$ be a free minimal Cantor $\Z^{2}$-system.
 For any $M \geq 1$, there
is a clopen set $Y \subset X$ such that the family
\[
P(x) = \{ m \in \Z^{2} \mid x \in \varphi^{m}(Y) \},
\]
for $x$ in $X$, 
is $M$-separated and $2M$-syndetic.
\end{prop}

\begin{proof}
For each $x$ in $X$, select a clopen set $V_{x}$ such that the sets
 $\varphi^{m}(V_{x}), m \in B(0, M)$ are pairwise disjoint. These sets form an open cover of $X$.
Select a finite set $ x_{1}, x_{2}, \ldots, x_{n} \in X$ such that 
$V_{x_{i}},  1 \leq i \leq n$ cover $X$. Let $Y_{1} = V_{x_{1}}$ and for each $i > 1$, let
\[
Y_{i} = Y_{i-1} \cup \left[ V_{x_{i}} \setminus ( \cup_{n \in B(0, M)} \varphi^{n}(Y_{i-1}) ) \right].
\]
Put $Y = Y_{n}$ and let $P(x), x \in X$, be as in the statement. It is easy to see that $P$ is $M$-separated. It is also easy to check that every point
in $\Z^{2}$ is distance at most $M$ from some point in $P(x)$. Since 
$M \geq 1$, it then follows that $P(x)$ is $2M$-syndetic, for any $x$ in $X$.
\end{proof}

\section{Refining tessellations}

In the last section, we gave a method of producing $\varphi$-regular tessellations.
The next step is to show how we may produce a sequence having larger and larger
elements (more and more separated) in such a way that each element of one is the union
of elements from the previous. At the same time, we will need several extra technical 
conditions which will be used  later in the proof of Theorem \ref{thm:affapp}. While we will provide
rigorous and fairly complete arguments, most of these properties can be seen fairly easily
by drawing some pictures.

Before stating the result, we will need some notation. This will also be used in later 
sections.
We are considering a tessellation, $\mathcal{T}$, of the plane by polygonal 
regions, with non-overlapping interiors.
Given a point $u$ in $\R^{2}$, we would like to say that this point belongs to a unique
element of $\mathcal{T}$. Of course, this is false since the elements overlap on their boundaries.
To resolve this difficulty in an arbitrary, but consistent way, we define, for any $u$
in $\R^{2}$ and polygon $t$, 
\[
u \in' t
\]
if, for all sufficiently small $\epsilon > 0$, we have 
\[
u + (\epsilon, \epsilon^{2}) \in t.
\]
For any subset $A \subset \R^{2}$, we define
\[
A \cap' t = \{ u \in A \mid u \in' t \}.
\]

A comment is in order regarding simplicial tessellations. In the process we are about to undertake, we
will take   unions of polygons, which may be disconnected. In addition, a vertex in some polygon may only belong to one other
element of the tessellation. So we would like to drop the terms `vertex' and `edge'.
 Instead, we 
would like to regard the `edges' in a combinatorial way as the (non-empty) intersection
of a pair of polygons. Geometrically, this set will be a union of line segments.
In a similar way, we would like to regard a (non-empty) three-way intersection
as a `vertex'. In fact, such a set may not be a single point.  Instead, we introduce the following notation.
For any simplicial tessellation $\mathcal{T}$ and  $1 \leq k \leq 3$, 
we let $\mathcal{T}^{k}$ denote the set of $k$-tuples, $(t_{1}, \ldots , t_{k})$ 
in $\mathcal{T}$ such that $\cap_{i=1}^{k} t_{i}$ is non-empty and $t_{i} \neq t_{j}$, 
for $i \neq j$.

\begin{thm}
\label{thm:5_refine}
Let $E \geq 3$ be as given in \ref{prop:3_tessellate}.
There exists a sequence of $\varphi$-regular simplicial tessellations $\mathcal{T}_{l}(x), x \in X, l \geq 0$,
satisfying each of the following conditions, for all $l \geq 0$ and $x$ in $X$:
\begin{enumerate}
\item
 $\mathcal{T}_{l}(x)$ has capacity $\max \{ l, E \}$, 
\item
 $\mathcal{T}_{l+1}(x)$ is $l +  diam(\mathcal{T}_{l})$-separated,
\item
each element of  $\mathcal{T}_{l}(x)$ meets at most $E$ (as in \ref{prop:3_tessellate}) other elements,
\item
if two elements of  $\mathcal{T}_{l}(x)$ meet, then there are exactly two others which meet both
of them,
\item
each element of  $\mathcal{T}_{l}(x)$ is contained in an element of  $\mathcal{T}_{l+1}(x)$,
\item
if $(t_{1}, t_{2}, t_{3})$ and $(t_{1}, t_{2}, t_{3}')$ are 
in $\mathcal{T}^{3}_{l+1}(x)$
with $t_{3} \neq t_{3}'$, 
letting \linebreak 
$\mathcal{N}^{3}(x,t_{1}, t_{2}, t_{3})$ denote the set of all $(t, t', t'')$ in
$\mathcal{T}_{l}^{3}(x)$ such that $t \subset t_{1}, t' \subset t_{2}, t'' \subset t_{3}$,
and $\mathcal{N}^{2}(x,t_{1}, t_{2}, t_{3})$ denote the set of all $t \subset t_{1}$ in 
$\mathcal{T}_{l}(x)$ such that: 
\begin{enumerate}
\item
there is $t' \subset t_{2}$ such that 
$(t,t')$ is in $\mathcal{T}_{l}^{2}(x)$, 
\item
for all $(t,t')$ in $\mathcal{T}_{l}^{2}(x)$, $t' \subset t_{1} \cup t_{2}$
and 
\item
$d(t, t_{3}) < d(t, t_{3}') - diam(\mathcal{T}_{l})$,
\end{enumerate}
we have 
\[
\# \mathcal{N}^{3}(x,t_{1}, t_{2}, t_{3}) < \# \mathcal{N}^{2}(x,t_{1}, t_{2}, t_{3}),
\]
\item
 for all  $t$ and $s$ in 
$\mathcal{T}_{l+1}(x)$, we have 
\[
2^{l} \# \{ n \in \Z^{2} \cap' t \mid d(n, \R^{2} \setminus t) \leq l \} \leq \# (\Z^{2} \cap' s).
\]
\end{enumerate}
\end{thm}

\begin{proof}
Let $b, a > 0$ be as in \ref{prop:3_tessellate}.
By \ref{prop:4_sep_synd}, we may choose a clopen, non-empty subset
of $X$ such that the associated $\varphi$-regular subset 
of $\R^{2}$, denoted $P_{0}(x), x \in X$,  is $E/b$-separated and $2E/b$-syndetic. 
We define $\mathcal{T}_{0}(x) = \mathcal{T}_{P_{0}(x)}$, for each $x$ in $X$, as in the last section. By \ref{prop:TP_locally_derived} and \ref{prop:4_local_derive}, it is $\varphi$-regular.
By \ref{prop:3_tessellate}, this has capacity $E$ and satisfies  properties 1, 3 and 4 (which are the only ones that do not involve $\mathcal{T}_{1}$).
We also obtain a map, denoted $P_{0}(x, t)$,
which to any $x$ in $X$ and element $t$ in $\mathcal{T}_{0}(x)$ assigns a point $u$ in
the interior of $t$.
This function is $\varphi$-regular in the following sense:
\begin{enumerate}
\item for $x$ in $X$ and $t$ in $\mathcal{T}_{0}(x)$, 
$P_{0}(\varphi^n(x),t+n)=P_{0}(x,t)+n$, 
\item for $x$ in $X$ and $t$ in $\mathcal{T}_{0}(x)$, 
there exists a clopen neighborhood $U$ of $x$ such that 
for any $y$ in $U$ we have $t$ in $\mathcal{T}_{0}(y)$,
 and $P_{0}(y,t)=P_{0}(x,t)$. 
\end{enumerate}

Next, we suppose that we have found a simplicial tessellation $\mathcal{T}_{l}(x), x \in X$
 satisfying the desired conditions, for some $l \geq 0$. Suppose also, that we have a point, 
$P_{l}(x,t)$, for $x \in X$, $t \in \mathcal{T}_{l}(x)$, which is in the interior of $t$. 
Moreover, this function is $\varphi$-regular in the sense above. Also suppose that, for any $x$ in $X$, 
the sets $\{ P_{l}(x, t) \mid t \in \mathcal{T}_{l}(x) \}$ are $l$-separated.
Let $D_{l}$ denote the diameter of $\mathcal{T}_{l}$.
Note that since each element of $\mathcal{T}_{l}(x)$ contains a ball of radius $l$, $D_{l} \geq 2l$.

Let $\rho = \frac{2D_{l}+1}{\sin(a/2)}$. Note that
$\rho \geq 2D_{l}+1$. Fix $x$ in $X$.
We consider, for any $u$ in $\R^{2}$,  the number of points in $B(u, 3\rho)$ of the form
$P_{l}(x,t)$. For such  points, the sets $B(P_{l}(x,t), l/2)$ 
are non-overlapping and are contained in $B(u, 3\rho + l/2)$. Hence, the number of
such points is bounded by some constant, which we denote by $K_{l} \geq 1$, depending on $l$, but not
$x$ or $u$. 
We find a constant $M$ satisfying each of the following:
\begin{eqnarray*}
bM  &  >  & 4\rho + 8K_{l}^{2}(D_{l} + 1) + l  \geq 8 D_{l} + l +1, \\
(bM)^{2} & > &  2^{l+2}(2(l + D_{l}) +2)E(8M + 2(l + D_{l}) + 2),
\end{eqnarray*}
noting that in the second inequality, the left hand side is quadratic in $M$, while the right hand side
is linear.

By Proposition \ref{prop:4_sep_synd}, we may find a clopen set $U \subset X$ such that the associated $P(x), x \in X$ is $2M$-syndetic and
$M$-separated. We let $\mathcal{T}_{P}(x), x \in X$ be the associated 
 simplicial tessellation and $P(x, \cdot)$ denote the canonical
map from $\mathcal{T}_{P}(x)$ to $\R^{2}$, as described in the section 3.
By Propositions \ref{prop:TP_locally_derived} and \ref{prop:4_local_derive}, $\mathcal{T}_{P}$ is $\varphi$-regular.

For each $x$ in $X$, $t$ in $\mathcal{T}_{P}(x)$, we define
\[
\tilde{t} = \bigcup_{t' \in \mathcal{T}_{l}(x), P_{l}(x,t') \in' t} t'.
\]
We now define a new  tessellation, denoted $\mathcal{T}_{l+1}(x)$ to be
the collection of all $\tilde{t}$, where $t$ is in  $\mathcal{T}_{P}(x)$. Since
$\mathcal{T}_{l}, \mathcal{T}_{P}$ are $\varphi$-regular and $P_{l}(\cdot, \cdot)$ is $\varphi$-regular in the sense described above,
this new tessellation is $\varphi$-regular. Clearly, condition 5 is satisfied. 

The first step is to observe that, for any $t$ in  $\mathcal{T}_{P}(x)$ as above,
each point in $\tilde{t}$ is in some element $t'$ of $\mathcal{T}_{l}(x)$ with
$P_{l}(x,t') \in' t$. As the diameter of $t'$ is at most $D_{l}$, it follows that
every point of $\tilde{t}$ is within distance $D_{l}$  of $t$.

We next verify that $\mathcal{T}_{l+1}(x)$ has capacity $\max \{ l+1, E \}$. First, each element
of $\mathcal{T}_{l+1}(x)$ contains an element of $\mathcal{T}_{l}(x)$ and hence has capacity
$E$. By Proposition \ref{prop:3_tessellate}, each $t$ in $\mathcal{T}_{P}(x)$
contains the ball $B(P(t), bM)$. As the diameter of $\mathcal{T}_{l}$ is $D_{l}$, the ball
$B(P(t), bM - 2D_{l})$ is contained in $\tilde{t}$. Since $bM - 2D_{l} \geq 6D_{l} + l + 1 \geq l+1$,
it follows that $\mathcal{T}_{l+1}$ has capacity $l+1$.

We know that any two disjoint elements, $t_{1}, t_{2}$ of $\mathcal{T}_{P}(x)$ are separated by
distance at least $bM$, by part 4 of Proposition \ref{prop:3_tessellate}. From the 
observation above, 
this means 
\[
d(\tilde{t_{1}}, \tilde{t_{2}}) \geq d(t_{1}, t_{2}) - 2D_{l} \geq bM - 2 D_{l} \geq 6D_{l} + l + 1 \geq l+D_{l}.
\]

We will show that the map sending $t$ in $\mathcal{T}_{P}(x)$ to $\tilde{t}$ in 
$\mathcal{T}_{l+1}(x)$ is a bijection which preserves non-trivial (multiple) intersections.
The first step in this is to observe from the last paragraph that if $t_{1}$ and $t_{2}$ are disjoint, then so are $\tilde{t_{1}}$ and $\tilde{t_{2}}$.

Now, we want to consider the situation that $t_{1}, t_{2}, t_{3}$ are three distinct
elements of $\mathcal{T}_{P}(x)$ with a non-trivial intersection, say at $u$. We will show that $\tilde{t_{1}}, \tilde{t_{2}}$ and $\tilde{t_{3}}$ also have non-trivial intersection. Consider the
closed disc, $D$, centred at $u$ with radius $\rho$ and let $C$ be the boundary of $D$. 
Consider any  $s$ 
in $\mathcal{T}_{P}(x)$, other than $t_{1}, t_{2}, t_{3}$. 
By part 6 of Proposition \ref{prop:3_tessellate}, the distance from any points of $D$ to $s$ is at least $bM - \rho$. It follows from $bM - \rho \geq D_{l} +1$ that $D$ does not meet $\tilde{s}$.

For each $ i=1, 2, 3$, let $u_{i}$ 
be the point on the circle in $t_{i}$ where the line from $u$ to $u_{i}$ bisects the angle
of $t_{i}$ at $u$. It follows from 5 and 6 of Proposition \ref{prop:3_tessellate},
the definition of $\rho$ and some simple trigonometry, that 
\[
d(u_{i}, t_{j}) \geq \rho \sin( \frac{a}{2}) = 2 D_{l} + 1,
\]
for $j \neq i$. In fact, this estimate holds replacing $u_{i}$ with any point on the arc
of $C$ between $u_{i}$ and $u_{i'}$ for $j \neq i, i'$.
This means that $u_{i}$ is in $\tilde{t_{i}}$ and no other
element of $\mathcal{T}_{l+1}(x)$. 
Consider the arc of $C$ from  $u_{1}$ to $u_{2}$. 
Arguments similar to those earlier show that this arc does not meet $\tilde{t_{3}}$, 
nor does it meet any elements of $\mathcal{T}_{l+1}(x)$ other than $\tilde{t_{1}}, \tilde{t_{2}}$.
The endpoints of the arc are contained in $\tilde{t_{1}}$ and  $\tilde{t_{2}}$, 
respectively. Since the arc is connected, there exists a point, $v$, on  the arc which is in both.
Thus, we find a point in $\tilde{t_{1}} \cap \tilde{t_{2}}$ and in no other element of
$\mathcal{T}_{l+1}(x)$. 

Finally, we want to show that $\tilde{t_{1}} \cap \tilde{t_{2}} \cap \tilde{t_{3}}$ is non-empty.
Let us assume that this intersection is empty. 
We consider all points inside of $C$, which lie in more than one of the sets. These are points
on the boundary of the $\tilde{t_{i}}$'s, which are the union of edges from the first
simplicial tessellation, $\mathcal{T}_{0}$, which is simplicial. Let $\mathcal{E}$ denote
these edges. In $\mathcal{T}_{0}$, at most three edges meet at each vertex. It is easy to see that
if our triple intersection is empty, no three edges of $\mathcal{E}$ can meet at a vertex. 
Similarly, any vertex which meets one edge of $\mathcal{E}$, must also meet another.
So each vertex of $\mathcal{T}_{0}$ inside $C$ must meet two elements of $\mathcal{E}$ or none. 
Consider all edges which meet $C$; 
these must be connected in pairs by  path in $\mathcal{E}$
passing through the inside of $C$. We conclude that there are an even number of such edges.

Now consider the arc of $C$ from $u_{1}$ to $u_{2}$. As noted above, it does not meet 
$\tilde{t}_{3}$, but begins in $\tilde{t}_{1}$ and ends in $\tilde{t}_{2}$.
It therefore meets an odd number of edges which separate 
 $\tilde{t}_{1}$ and  $\tilde{t}_{2}$. In a similar way, the other two arcs of $C$ also contain
an odd number of edges of $\mathcal{E}$. But this means that there are an odd number of such edges.
This contradiction establishes the desired result.

We are now ready to show that the map sending $t$ in $\mathcal{T}_{P}(x)$ to $\tilde{t}$ in $\mathcal{T}_{l+1}(x)$ preserves non-trivial intersections. If $t_{1}$ and $t_{2}$ have a non-trivial intersection, then there is a $t_{3}$ which meets both of them. The argument above shows that $\tilde{t_{1}}, \tilde{t_{2}}$ and $\tilde{t_{3}}$ have non-trivial intersection and so $\tilde{t_{1}}$ and $\tilde{t_{2}}$ do also. We have already established the converse. We have also established the three-way intersection version of this. For the converse for three-way intersections, if $\tilde{t_{1}}, \tilde{t_{2}}$ and $\tilde{t_{3}}$ have non-trivial intersection, then each pair has non-trivial intersection. Then the same holds for each pair of $t_{1}, t_{2}$ and $t_{3}$. Then since $\mathcal{T}_{P}(x)$ is simplicial, $t_{1}, t_{2}$ and $t_{3}$ have non-trivial intersection.  This same argument shows that $\mathcal{T}_{l+1}(x)$ has no non-trivial four-way intersections. This completes 
 the proof.

It follows easily from this that $\mathcal{T}_{l+1}(x)$ is simplicial and that conditions 3 and 4 hold. Moreover, if $\tilde{t_{1}}$ and $\tilde{t_{2}}$ are two elements of $\mathcal{T}_{l+1}(x)$ with trivial intersection, it follows that $t_{1}$ and $t_{2}$ also have trivial intersection. We have then established above that $d(\tilde{t_{1}}, \tilde{t_{2}}) \geq  l+D_{l}$ and so
 $\mathcal{T}_{l+1}(x)$ is $l + diam(\mathcal{T}_{l})$-separated.

We next consider condition 6. Begin with $t_{1},t_{2},t_{3},t_{3}'$ 
in $\mathcal{T}_P(x)$, 
so that $\tilde{t_{1}},\tilde{t_{2}},\tilde{t_{3}},\tilde{t_{3}'}$ 
are as in the hypothesis of 6. 
Let $D$ be the closed disk centred at $t_{1}\cap t_{2}\cap t_{3}$ 
with radius $\rho+2K_{l}^{2}(D_{l}+1)$. 
By part 6 of Proposition \ref{prop:3_tessellate}, $D$ is contained in $t_{1} \cup t_{2} \cup t_{3}$ and $D \cap t_{i}$ is a sector, for $i = 1, 2, 3$. By the same argument as before, we see that $D$ does not meet any element of $\mathcal{T}_{l+1}(x)$ other than $\tilde{t_{1}}, \tilde{t_{2}}$ and $\tilde{t_{3}}$.
Consider $D \cap t_{1}\cap t_{2}$, 
which is a line segment of length $\rho+2K_{l}^{2}(D_{l}+1)$. 
Consider the interval 
which begins at distance  $\rho$ from $t_{1} \cap t_{2} \cap t_{3}$ 
and ends at the boundary of $D$. 
Its length is clearly $2K_{l}^{2}(D_{l}+1)$. 
Hence, for each $k=1,2,\dots,2K_{l}^{2}$, 
we may place a line segment, denoted $L_{k}$, with length $2D_{l}$, 
perpendicular to $t_{1} \cap t_{2}$, midpoint on $t_{1} \cap t_{2}$ and 
so that any two are distance at least $D_{l}+1$ apart. 
It is an easy geometric argument to see that 
for $u$ in any of these, we have 
\[ 
d(u,t_{3}) < \rho + 2K_{l}^{2}(D_{l}+1) < \frac{bM}{4}. 
\]
Since $d(t_{3},t_{3}')\geq bM$, we also get $d(u,t_{3}')\geq \frac{3bM}{4}$. 
Therefore, we have 
\begin{eqnarray*} 
d(u,\tilde{t_{3}})  & <  & \frac{bM}{4}+D_{l},  \\
d(u,\tilde{t'_{3}}) & \geq &  \frac{3bM}{4}-D_{l}. 
\end{eqnarray*}
An endpoint of $L_{k}$ is in one of $t_{1},t_{2}$ and 
is distance $D_{l}$ from the other. It follows that one endpoint of $L_{k}$ is contained in $\tilde{t}_{1}$ and the other is contained in $\tilde{t}_{2}$. Since $L_{k}$ is connected, there must be some point on $L_{k}$ which is contained in both $\tilde{t}_{1}$ and $\tilde{t}_{2}$.  This point is then contained in some pair $(t,t')$ in $\mathcal{T}_{l}^{2}(x)$ with $t \subset \tilde{t}_{1}$ and $t' \subset \tilde{t}_{2}$. These are distinct for different $k$
because the $L_{k}$ are separated by distance at least $D_{l}+1$. 
Besides, from $\frac{bM}{2} > 4D_{l}$, we have 
\[ 
d(t,\tilde{t_{3}}) <  \frac{bM}{4} + D_{l}
  <  \frac{3bM}{4} - 3D_{l} \leq d(t,\tilde{t'_{3}}) - D_{l}. 
\]
The remaining properties of these pairs are easily checked from the conditions given above.
From this we conclude that the second number in condition 6 is at least $2K_{l}^{2}$. 
Now we show that the first  number
is at most $2K_{l}(K_{l}-1)$. If $t, t', t''$ are as given, then $t'' \subset \tilde{t_{3}}$ implies that $P_{l}(x,t'')$ 
is in $t_{3}$. It follows that, letting $u$ be any point of  $ t \cap t''$,  
\[
d(P_{l}(x,t), t_{3})  \leq    d(P_{l}(x,t),u) + d(u, P_{l}(x,t'')) \leq
2D_{l} . 
\]
Similarly, we have 
\[
 d( P_{l}(x, t'), t_{3}) \leq 2D_{l}.
\]
In addition, we have
\begin{eqnarray*}
P_{l}(x, t) \in t_{1}, & &  P_{l}(x, t') \in t_{2}, \\
d(P_{l}(x, t), t_{2}) & \leq & 2D_{l}, \\
d(P_{l}(x, t'), t_{1}) &  \leq  &  2D_{l}.
\end{eqnarray*}
Finally, it is an easy geometric exercise, using the fact that the angles made at
 $\{ u \} = t_{1} \cap t_{2} \cap t_{3}$ are at least $a$, to see that these conditions imply 
\[
P_{l}(x,t), P_{l}(x, t') \in B(u, 3\rho).
\]
The conclusion follows at once from the definition of $K_{l}$ and the fact that for a given $t,t'$
there are at most two $t''$ which meet both $t$ and $t'$.

We now consider condition 7. Let $t$ be in $\mathcal{T}_{P}(x)$ and show the
condition holds for $\tilde{t}$. We know  
that $t$ is contained in $B(P(x,t), 4M)$ and hence by 
 Proposition \ref{prop:3_tessellate} has at most $E$ edges, each of which
is of length at most $8M$. If $n$ is in $\Z^{2} \cap' \tilde{t}$ and $d(n, \R^{2} \setminus \tilde{t}) \leq l$, then 
$d(n, u) \leq l + D_{l}$, for some point $u$ on an edge of $t$. We cover each edge of $t$ 
by a rectangle of width $2(l + D_{l})$ and length $\lambda + 2(l + D_{l})$, where $\lambda$ is the
length of the edge. These rectangles will contain all points $n$ which we consider above.
We need to estimate the number of lattice points contained in a rectangle in the plane. At each such point, place a unit square, centred at the points. As these are non-overlapping, the area of their union is equal to the number of points. Moreover, each is contained in a rectangle whose side lengths are 2 larger than the original.
Hence, we see that the number of lattice points covered by a rectangle is bounded above
by the product of its length plus 2 and its width plus 2. Therefore, we have 
\begin{eqnarray*}
\# \{ n \in \Z^{2} \cap' \tilde{t} \mid d(n, \R^{2} \setminus \tilde{t})  \leq l \} & 
          \leq  &        (2(l + D_{l}) + 2)E(8M + 2(l + D_{l}) + 2) \\
          & \leq &  \frac{ (bM)^{2}}{2^{l+2}},
\end{eqnarray*}
from our choice of $M$. Since $\mathcal{T}_{l+1}$ has capacity 
\[
 bM - 2D_{l} = \frac{bM}{2} + \frac{bM}{2} - 2D_{l} \geq \frac{bM}{2},
\]
 $ s$ contains
a ball of radius $\frac{bM}{2}$, which, in turn, contains a square of side length $\frac{bM}{2} + 1$ and hence,
it contains at least $\frac{(bM)^{2}}{4}$ points of $\Z^{2}$ in its interior. Thus, we have
\[
2^{l} \# \{ n \in \Z^{2} \cap' \tilde{t} \mid d(n, \R^{2} \setminus \tilde{t}) \leq l \} \leq \frac{(bM)^{2}}{4} \leq \# (\Z^{2} \cap' s).
\]

As a final point, we define the function $P_{l+1}(x, \tilde{t})= P(x, t)$, for all
$x$ in $X$, $t$ in $\mathcal{T}_{P}(x)$. Notice that by part 3 of Proposition \ref{prop:3_tessellate},
 $B(P(x,t),bM)$ is contained in $t$. Since $bM > D_{l}+1$, $P(x,t)$ is in the interior 
of $\tilde{t}$ as required.

\end{proof}

\section{AF-equivalence relations and boundaries}

Here, we want to use our earlier construction of a nested sequence of $\varphi$-regular tessellations to give a proof of Theorem \ref{thm:affapp}. 
This needs, first of all, an AF-relation. The obvious choice is by using the interiors of the cells in the tessellation. These equivalence relations are actually too large. We will refine them. The interior of a tile $t$ will be subdivided by considering all $t', t''$ such that $t, t', t''$ have non-trivial intersection. For each such triple, we will determine a subset of the lattice points in $t$. Although the construction is completely combinatorial, it is reasonable to imagine these sets geometrically in the following way:

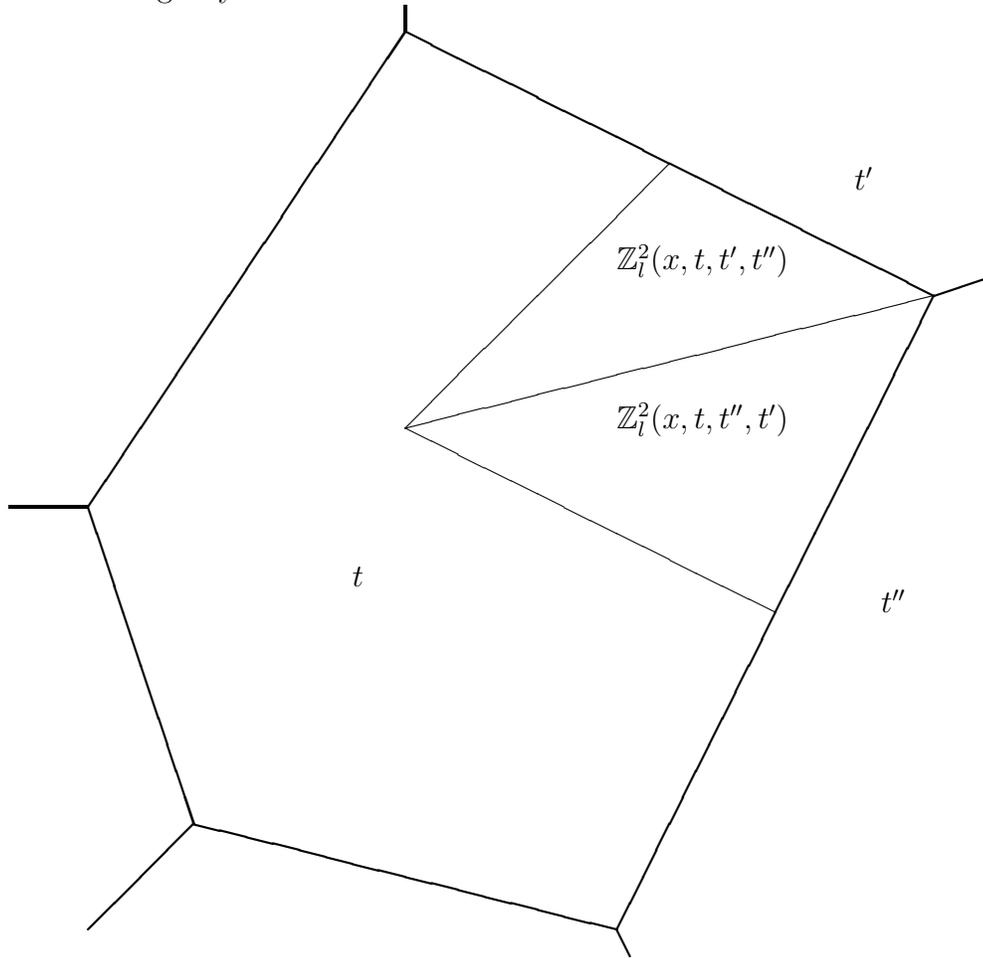
\begin{picture}(370,360)(0,0)

\thicklines

\put(0,170){\line(1,0){30}}
\put(30,170){\line(2,3){120}}
\put(150,350){\line(0,1){10}}
\put(150,350){\line(2,-1){200}}
\put(350,250){\line(3,1){20}}
\put(230,10){\line(1,2){120}}
\put(230,10){\line(1,-2){5}}
\put(70,50){\line(4,-1){160}}
\put(30,10){\line(1,1){40}}
\put(30,170){\line(1,-3){40}}

\thinlines

\put(150, 200){\line(1,1){100}}
\put(150, 200){\line(4,1){200}}
\put(150, 200){\line(2,-1){140}}

\put(130,140){$t$}
\put(320,290){$t'$}
\put(330,130){$t''$}

\put(230,260){$\Z^{2}_{l}(x,t,t',t'')$}
\put(230,200){$\Z^{2}_{l}(x,t,t'',t')$}

\end{picture}

At the same time, we also keep track of boundary sets, $B^{2}$ along the edges and $B^{3}$ around the vertices.
The first application of the absorption theorem enlarges the equivalence relation along the edge boundaries, the second around the vertex boundaries.
Although we have not checked all the details, it seems likely that all of these constructions can be extended to the case of $\Z^{d}$-actions for $d > 2$. The single missing ingredient is the analogue of Proposition \ref{prop:3_tessellate} (and the Lemmas leading to it).

We begin with our refining sequence of $\varphi$-regular simplicial 
tessellations $\mathcal{T}_{l}, l \geq 0$ provided by Theorem \ref{thm:5_refine}  of the last section.
For any $x$ in $X$, we let $i(x, \cdot): \mathcal{T}_{l}(x) \rightarrow
\mathcal{T}_{l+1}(x)$ be the unique function such that $i(x,t) \supset t$, for every $t$ in 
$\mathcal{T}_{l}(x)$. If $k \geq 1$, we let $i^{k}$ denote the composition of
$k$ functions $i$ mapping $\mathcal{T}_{l}$ to $\mathcal{T}_{l+k}$, for any $l$.

For each $x$ in $X$, the tessellation $\mathcal{T}_{l}(x)$ naturally partitions (via $\in'$) the integer lattice into
finite equivalence relations. Each of these classes is indexed by an element of $\mathcal{T}_{l}(x)$.
However, we need to refine this equivalence relation.
In our refined relation, each equivalence class will be indexed by  elements 
$(t,t',t'')$ from $\mathcal{T}_{l}^{3}(x)$.
The union over fixed $t', t''$ will yield the elements of $ \Z^{2}  \cap' t$. Some care must be taken so
that this is done in a $\varphi$-regular way.
Let  $\mathcal{T}(x), x \in X$ be a $\varphi$-regular tessellation. For $x$ in $X$ and
any $t$ in $\mathcal{T}(x)$, we let $N(x, t)$ denote the set of all $t'$ in 
$\mathcal{T}(x)$, including $t$, which intersect $t$. 
 Consider all possible 
$(t, N(x,t))$, where $x$ is in $X$ and 
$t$ is in
 $\mathcal{T}_{0}(x)$. We consider $(t, N(x,t))$ and
$(t', N(x',t'))$  to be equivalent if they are translates of one another,
namely that there is a $u$ in $\R^{2}$ such that $t' = t + u$ and
$N(x',t') = N(x,t) +u$.
Since $\mathcal{T}_{0}$ is $\varphi$-regular, there are a finite number of equivalence classes.
We let $\mathcal{P}$ be a finite set containing exactly one representative of each
equivalence class.

Let $(t,N)$ be in  $\mathcal{P}$. For each $t' \neq t''$ in $N \setminus \{ t \}$ with $t \cap t' \cap t''$ non-empty , 
we define $b(t,N, t', t'')$ to be any point in $\Z^{2} \cap'  t $.
These should be chosen to be distinct for different ordered pairs $t', t''$. 
To see this is possible, note that
for any $x$, $\mathcal{T}_{0}(x)$ has capacity $E$, so each element
contains a ball of radius $E \geq 3$, hence a square of side length $E+1$ and so $E^{2} \geq 2E$ points
of the integer lattice.   The result then follows from condition 3 of Theorem \ref{thm:5_refine}.
Next, we partition
the elements of $ \Z^{2}  \cap' t$ into sets, $\Z^{2}_{0}(t, N, t', t'')$, indexed by the pairs,
 $t', t''$ as above. These should be chosen so that
$\Z^{2}_{0}(t, N, t', t'')$ contains $b(t,N,t',t'')$, for all $t', t''$.

Having chosen these items for our representative patterns $\mathcal{P}$, we extend their 
definition by translation as follows.
Let $x$ be in $X$ and $(t, t', t'')$ be in $\mathcal{T}_{0}^{3}(x)$. We find the unique
$n$ in $ \Z^{2}$ with $(t -n, N(x,t)-n)$ in $\mathcal{P}$ and define
\begin{eqnarray*}
b(x, t, t', t'') & =  & b( t-n, N(x,t)-n, t' -n, t'' -n) + n, \\
\Z^{2}_{0}(x,t,t', t'') & = & \Z^{2}_{0}(t-n, N(x,t)-n, t' -n, t'' -n) + n.
\end{eqnarray*}
Further, we define, for each $x$  in $X$ and $(t, t', t'')$ in $\mathcal{T}_{0}^{3}(x)$,
\[
B_{0}^{3}(x,t,t',t'') = \{ b(x,t,t',t'') \},
\]
and for  each $x$  in $X$ and $(t, t')$ in $\mathcal{T}_{0}^{2}(x)$,
\[
B_{0}^{2}(x, t, t') = \{ b(x, t, t', t'') \mid \text{ there is } t'' \text{ such that }
                        (t,t',t'') \in \mathcal{T}_{0}^{3}(x) \}.
\]
We also define analogues of these sets in $X$ by
\[
B_{0}^{2} = B_{0}^{3} = \{ \varphi^{b(x,t,t',t'')}(x) \mid x \in X, (t,t',t'') \in \mathcal{T}_{0}^{3}(x) \}.
\]

Although they have no specific geometric property, we refer to these points, $b(x, t, t', t'')$, as
`boundary points'. The following Lemma follows at once from the definitions and we omit the proof.

\begin{lemma}
For any $x$ in $X$, $(t, t', t'')$  in $\mathcal{T}_{0}^{3}(x)$ and $n$ in $\Z^{2}$, we have 
\begin{eqnarray*}
b(\varphi^{n}(x), t+n, t'+n, t''+n) & = & b(x, t, t', t'') +n, \\
B_{0}^{3}(\varphi^{n}(x), t+n, t'+n, t''+n) & = & B_{0}^{3}(x, t, t', t'') +n, \\
B_{0}^{2}(\varphi^{n}(x), t+n, t'+n) & = & B_{0}^{2}(x, t, t') +n, \\
\Z^{2}_{0}(\varphi^{n}(x), t+n, t'+n, t''+n) & = & \Z^{2}_{0}(x, t, t', t'') +n.
\end{eqnarray*}
\end{lemma}

Notice that $B_{0}^{2}$ and $B_{0}^{3}$ are clopen because of $\varphi$-regularity and Proposition \ref{prop:4_regular}.

Ultimately, we will need to enlarge our AF-equivalence relation by including
equivalences between boundary points. We can define these relations as follows.
We define
\begin{eqnarray*}
K_{0}^{3} & = & \{ (\varphi^{b(x,t,t',t'')}(x), \varphi^{b(x,t'',t,t')}(x)), \\
          &   &         (\varphi^{b(x,t'',t,t')}(x), \varphi^{b(x,t,t',t'')}(x)), \\
          &  &            (\varphi^{b(x,t,t',t'')}(x), \varphi^{b(x,t,t',t'')}(x)) \mid  \\
          &  &   x \in X,  (t, t', t'') \in \mathcal{T}_{0}^{3}(x) \}, \\
K_{0}^{2} & = & \{ (\varphi^{b(x,t,t',t'')}(x), \varphi^{b(x,t',t,t'')}(x)), \\ 
          &  &   (\varphi^{b(x,t,t',t'')}(x), \varphi^{b(x,t,t',t'')}(x))\mid   \\
          &  &  x \in X,  (t, t', t'') \in \mathcal{T}_{0}^{3}(x) \}.
\end{eqnarray*}

\begin{lemma}
\label{lemma:6_K_compact}
For $i = 2, 3$, $K_{0}^{i}$ is a compact open subequivalence relation of
$R_{\varphi} \mid B_{0}^{i}$.
\end{lemma}

\begin{proof}
That $K_{0}^{i}$ is an equivalence relation follows from the fact that if $(t, t', t'')$
is an element of $\mathcal{T}_{0}^{3}(x)$, then so is any permutation of these three
elements of $\mathcal{T}_{0}(x)$. The rest of the proof follows from the continuity of
$b(x,t,t',t'')$, which can be stated as follows. For any $x$, since $\mathcal{T}_{0}$ is
$\varphi$-regular, for any $(t, t', t'')$ in $\mathcal{T}_{0}^{3}(x)$,
there is a clopen neighbourhood  $U$ of $x$ such that
$(t, t', t'')$ is in $\mathcal{T}_{0}^{3}(x')$, for all $x'$ in $U$. Then from the
definition of $b$, we have that $b(x',t,t',t'') = b(x, t, t', t'')$, for all such $x'$.
It follows at once that $K_{0}^{i}$ is open. By allowing $x$ to vary over $X$, the sets $U$ obtained form 
an open cover. By selecting a finite subcover, it follows quite easily that $K_{0}^{i}$ is also compact.
\end{proof}

Next, we assume that we have defined, for some $l \geq 0$,
sets $\Z^{2}_{l}(x, t, t', t'')$ and $B^{3}_{l}(x,t,t', t'')$, for $x$ in $X$,
 $(t, t', t'')$ in $\mathcal{T}_{l}^{3}(x)$, and $B^{2}_{l}(x,t,t')$, 
for $x$ in $X$,
 $(t, t')$ in $\mathcal{T}_{l}^{2}(x)$. 

We have a simplicial tessellation 
$\mathcal{T}_{l+1}$ which refines $\mathcal{T}_{l}$.
The next step in our construction involves determining the sets 
$\Z^{2}_{l+1}(x,t,t', t'')$.

For each $x$ in $X$, we define a function 
\[
i_{3}(x, \cdot): \mathcal{T}_{l}^{3}(x) \rightarrow \mathcal{T}_{l+1}^{3}(x),
\]
and then set  
\[
\Z_{l+1}^{2}(x, s, s', s'') = \bigcup_{i_{3}(x,t,t',t'')=(s,s',s'')} \Z_{l}^{2}(x, t, t', t'').
\]
To define $i_{3}$, we consider the sets $\{ i(x,t), i(x,t'), i(x,t'') \}$ and 
$i(x,N(x,t))$. Notice that the
first is clearly contained in the second and both contain $i(x,t)$.
For $x$ in $X$ and $t$ in $\mathcal{T}_{l}(x)$, the facts that $\mathcal{T}_{l+1}(x)$ is
$diam( \mathcal{T}_{l})$-separated, by condition 2 of  Theorem \ref{thm:5_refine}, and is simplicial means that, for any $t$ in 
$\mathcal{T}_{l}(x)$, $i(x, N(x,t))$ consists of at most three elements of $\mathcal{T}_{l+1}(x)$. We first consider the case that $i(x, N(x,t))$ contains only
$i(x,t)$. Then we define $i_{3}(x,t, t', t'') = (i(x,t), s', s'')$, where $s', s''$ are any elements of $\mathcal{T}_{l+1}(x)$ such that $(i(x,t), s', s'')$ is in $\mathcal{T}_{l+1}^{3}(x)$.
Of course, this must be done in a $\varphi$-regular fashion; that is the choice depends only on the pattern in $\mathcal{T}_{l+1}(x)$ around $i(x,t)$, in particular, it depends only on $i(x,t)$ and $N(x, i(x,t))$.
Next, consider the case $i(x, N(x,t)) = \{ i(x,t), s' \}$, for some $s' \neq i(x,t)$. There are two  $s''$
such that $(i(x,t),s',s'')$ is in $\mathcal{T}_{l+1}^{3}(x)$ and we let $s''$ denote the one closest to $t$. In the case that they are equidistant from $t$, either may be chosen, but it should be done in a local way. We then set $i_{3}(x,t, t', t'') = (i(x,t), s', s'')$. 
Next, we consider the case that $i(x,t) = i(x, t') = i(x,t'')$, but $i(x, N(x,t)) = \{ i(x,t), s', s'' \}$, for some $s', s''$. It readily follows that $(i(x,t), s', s'')$ is in $\mathcal{T}_{l+1}^{3}(x)$. We then set $i_{3}(x,t,t', t'') = (i(x,t), s', s'')$, in the order so that they appear in clockwise fashion. (Any other local rule would work as well.) Next, we consider the case $\{ i(x,t), i(x,t'), i(x, t'') \} = \{ i(x,t), s' \}$, while $i(x, N(x,t)) = \{ i(x,t), s', s'' \}$. In this case, we set $i_{3}(x,t,t', t'') = (i(x,t), s', s'')$. Finally, we are left to consider the case that $i(x,t), i(x,t')$ and $i(x,t'')$ are all distinct. In this case, we set $i_{3}(x,t,t',t'') = (i(x,t),i(x,t'),i(x,t''))$.

We first establish the following. Its proof is an easy consequence of the definition and we omit it.

\begin{lemma}
\label{lemma:6_i3}
\begin{enumerate}
\item If $(t,t',t'')$ is in $\mathcal{T}_{l}^{3}(x)$ and 
$i(x,t),i(x,t'),i(x,t'')$ are all distinct, then 
$i_{3}(x,t,t',t'')=(i(x,t),i(x,t'),i(x,t''))$. 
\item If $(t,t',t'')$ is in $\mathcal{T}_{l}^{3}(x)$ and 
$i(x,t)\neq i(x,t')$, then \linebreak 
$i_{3}(x,t,t',t'')=(i(x,t),i(x,t'),s)$ for some $s$ in $\mathcal{T}_{l+1}(x)$. 
\item If $(t,t',t_{i})$ are in $\mathcal{T}_{l}^{3}(x)$ for $i=1,2$ and 
$i(x,t)\neq i(x,t')$, then \linebreak
$i_{3}(x,t,t',t_{1})=i_3(x,t,t',t_{2})$. 
\end{enumerate}
\end{lemma}

For $(t,t',t'')$ in $\mathcal{T}^{3}_{l+1}(x)$ or $(t, t')$ in $\mathcal{T}^{2}_{l+1}(x)$, we define $B^{3}_{l+1}(x,t,t',t'')$ and $B^{2}_{l+1}(x, t, t')$, respectively, by
\begin{eqnarray*}
B^{3}_{l+1}(x,t,t',t'') & = & 
 \bigcup_{i(x,s)=t, i(x,s')=t', i(x,s'')=t''}  
             B^{3}_{l}(x, s, s', s''), \\
B^{2}_{l+1}(x, t, t') & = & 
 \bigcup_{i(x,s) = t, i(x,s') = t'} B^{2}_{l}(x,s,s'). 
\end{eqnarray*}
We transfer these sets to $X$ as follows. For each $l \geq 0$, let
\begin{eqnarray*}
B_{l}^{3} & = & \{ \varphi^{b}(x) \mid b \in B^{3}_{l}(x, t, t', t''), 
                          (t, t', t'') \in \mathcal{T}_{l}^{3}(x) \}, \\
B_{l}^{2} & = & \{ \varphi^{b}(x) \mid b \in B^{2}_{l}(x, t, t'), 
                          (t, t') \in \mathcal{T}_{l}^{2}(x) \}.
\end{eqnarray*}

\begin{lemma}
\label{lemma:6_B}
Let $x$ be in $X$ and let $(t, t', t'')$ be in $\mathcal{T}_{0}^{3}(x)$. 
\begin{enumerate}
\item
For any 
$l \geq 1$, $\varphi^{b(x,t,t',t'')}(x) $ is in $B^{3}_{l}$ if and only if 
\linebreak
$i^{l}(x,t), i^{l}(x,t'), i^{l}(x,t'')$ are all distinct. Moreover, in this
case, we have 
\[
b(x,t,t',t'') \in B_{l}^{3}(x, i^{l}(x,t), i^{l}(x,t'), i^{l}(x,t'')).
\]
\item
For any 
$l \geq 1$, $\varphi^{b(x,t,t', t'')}(x)  $ is in $B^{2}_{l}$ if and only if 
$i^{l}(x,t), i^{l}(x,t')$ are  distinct. Moreover, in this
case, we have 
\[
b(x,t,t',t'') \in B_{l}^{2}(x, i(x,t), i(x,t')).
\]
\item
The sets $B^{3}_{l}$ and $B^{2}_{l}$ are clopen. We also have $B^{3}_{l} \subset B^{2}_{l}$, 
$B^{3}_{l+1} \subset B^{3}_{l}$ and $B^{2}_{l+1} \subset B^{2}_{l}$. 
\item
For $i = 2,3$, $B_{l}^{i}$ 
is invariant under $K_{0}^{i}$.
\end{enumerate} 
\end{lemma}

\begin{proof}
We prove the first statement by induction on $ l \geq 0$ (we regard $i^{0}$ as the identity
map). It is clearly valid for $l=0$. Now assume it is true for $l$ and suppose that
$\varphi^{b(x,t,t',t'')}(x)$ is in $B_{l+1}^{3}$. This means that $b(x,t,t',t'')$ is in
$B_{l+1}^{3}(x,s,s',s'')$ for some $(s,s',s'')$ in $\mathcal{T}_{l+1}^{3}(x)$. Then by  definition,
$b(x,t,t',t'')$ is in $B_{l}^{3}(x,u,u',u'')$ for some $(u,u',u'')$ 
in $\mathcal{T}_{l}^{3}(x)$ with $i(x,u)=s, i(x,u')=s', i(x,u'')=s''$. 
It follows from the induction hypothesis
that $i^{l}(x,t), i^{l}(x,t'), i^{l}(x,t'')$ are all distinct and $b(x,t,t',t'')$ is in \linebreak
$B_{l}^{3}(x, i^{l}(x,t), i^{l}(x,t'), i^{l}(x,t''))$. But since the sets $B^{3}_{l}(x,w,w',w'')$ are disjoint
for different $(w, w', w'')$, it follows that $(i^{l}(x,t), i^{l}(x,t'), i^{l}(x,t'')) = (u,u',u'')$. Then  we
have
\begin{eqnarray*}
(i^{l+1}(x,t), i^{l+1}(x,t'), i^{l+1}(x,t''))  
& = & (i(x,u),i(x,u'),i(x,u'')) \\ 
  & =  & (s, s', s'')
\end{eqnarray*}
and this completes the proof. For the converse direction, if \linebreak $i^{l+1}(x,t), i^{l+1}(x,t'), i^{l+1}(x,t'')$ 
are all distinct, then clearly, so are $i^{l}(x,t)$, \linebreak
 $i^{l}(x,t'), i^{l}(x,t'')$. Then we may apply the induction 
hypothesis to conclude that $b(x,t,t',t'')$ is in $B_{l}^{3}(x,i^{l}(x,t), i^{l}(x,t'), i^{l}(x,t'') )$.
It then follows from the definition and the fact that $i^{l+1}(x,t), i^{l+1}(x,t'), i^{l+1}(x,t'')$ are distinct
that $b(x,t,t',t'')$ is in $B_{l+1}^{3}(x, i^{l+1}(x,t), i^{l+1}(x,t'), i^{l+1}(x,t''))$. 
The conclusion follows from the definition of $B_{l+1}^{3}$.

The proof of the second statement is completely analogous to the first and
 we omit it. The third statement is clear from the definitions. The last part is clear from the definition of
$K_{0}^{i}$ and the first two parts which we have already established.
\end{proof}

We define  subequivalence relations, $R_{l} \subset R_{\varphi}$, for $l \geq 0$,  by
\[
R_{l}=  
\{ (\varphi^{m}(x), \varphi^{n}(x)) \mid x \in X, (t,t',t'') \in \mathcal{T}_{l}^{3}(x),
   m,n \in \Z^{2}_{l}(x, t, t', t'') \}.
\]

\begin{prop}
\label{prop:6_R}
For each $l \geq 0$, $R_{l}$ is a  compact, open subequivalence relation 
of $R_{\varphi}$. Moreover, we have $R_{l}\subset R_{l+1}$.
\end{prop}

\begin{proof}
 It is clear that $R_{l}$ is a subequivalence relation of $R_{\varphi}$. Since we defined $\Z^{2}_{l}(x, t, t', t'')$ in a $\varphi$-regular fashion, $R_{l}$ is open. Compactness of $R_{l}$ follows from the fact that $\Z^{2}_{l}(x,t,t',t'')$ is a subset of $\Z^{2}$ whose diameter is bounded, uniformly over all $x$ and $(t,t',t'')$. Finally, since $\Z^{2}_{l+1}(x,s, s', s'')$ is a union of sets of the form $\Z^{2}_{l}(x, t, t', t'')$, we have $R_{l} \subset R_{l+1}$.
\end{proof}

We define
\[
R = \cup_{l \geq 0} R_{l}.
\]
By \ref{prop:6_R}, $R$ is an open AF-subequivalence relation of $R_{\varphi}$.

We also define
\begin{eqnarray*}
B^{3}  &  = &  \cap_{l \geq 0} B^{3}_{l}, \\
B^{2} & = & \cap_{l \geq 0} B^{2}_{l}.
\end{eqnarray*}
It is clear that $B^{3} \subset B^{2}$.
Moreover, we define the equivalence relations $K^{3}$ and $K^{2}$ to be the restrictions of
$K^{3}_{0}$ and $K^{2}_{0}$ to the sets $B^{3}, B^{2}$, respectively.

\begin{prop}
\label{prop:6_B2_etale}
$B^{2}$ is an \'{e}tale subset of $X$ for the relation $R$.
\end{prop} 

\begin{proof} 
We suppose that $(x_{1}, x_{2})$
is in $R$, with $x_{1}, x_{2}$ in $B^{2}$. We wish to find an open subset $U$ of $R$
where the maps $r, s$ are local homeomorphisms and such that, for any $(x_{1}', x_{2}')$ in $U$,
$x'_{1}$ is in $B^{2}$ if and only if $x'_{2}$ is also in $B^{2}$. First of all, we may find $L \geq 1$
such that $(x_{1}, x_{2})$ is in $R_{L}$. 
Since $x_{1}$ and $x_{2}$ are in $B^{2} \subset B^{2}_{0}$ and since they are equivalent in
$R \subset R_{\varphi}$, we may find $x$ in $X$, $(t_{1}, t_{1}', t_{1}''), (t_{2}, t_{2}', t_{2}'')$
in $\mathcal{T}^{3}_{0}(x)$ such that
\[
x_{1} = \varphi^{b(x, t_{1}, t_{1}', t_{1}'')}(x), x_{2} = \varphi^{b(x, t_{2}, t_{2}', t_{2}'')}(x).
\]
We also choose a clopen neighbourhood, $U_{1}$, of $x$ sufficiently small so that 
\[
U =  \{ (\varphi^{b(x, t_{1}, t_{1}', t_{1}'')}(x'), \varphi^{b(x, t_{2}, t_{2}', t_{2}'')}(x'))
\mid x' \in U_{1} \}
\]
is contained in $R_{L}$.

 From the regularity of the tessellations $\mathcal{T}_{l}$, $ 0 \leq l \leq L$, we may also assume that the neighbourhood $U_{1}$ of $x$ is chosen sufficiently small so that, for any $x'$ in $U_{1}$, 
$(t_{1}, t_{1}', t_{1}''), (t_{2}, t_{2}', t_{2}'')$ are in $\mathcal{T}_{0}^{3}(x')$ and
\begin{eqnarray*}
b(x,t_{j},t'_{j}, t''_{j}) &  =  & b(x',t_{j},t'_{j}, t''_{j}), \\
i^{l}(x, t_{j}) & = & i^{l}(x', t_{j}), \\
 i^{l}(x, t_{j}') &  = &  i^{l}(x', t_{j}'),\\
 i^{l}(x, t_{j}'') &  = &  i^{l}(x', t_{j}''),
\end{eqnarray*}
for all $ 1 \leq l \leq L, j = 1,2$.

Since $(x_{1}, x_{2})$ are in $R_{L}$, we have $b(x, t_{1}, t_{1}', t_{1}'')$ and
$b(x, t_{2}, t_{2}', t_{2}'')$ are in $\Z^{2}_{L}(x, s, s', s'')$ for some  $(s, s', s'')$
in $\mathcal{T}_{L}^{3}(x)$. By Lemma \ref{lemma:6_B}, since $x_{1}$ and $x_{2}$ are in $B^{2}$, we know
that
$i^{l}(x, t_{1}) \neq i^{l}(x, t_{1}'), i^{l}(x, t_{2}) \neq i^{l}(x, t_{2}')$ for all $l \geq 1$.
In particular, we have 
\[
i^{L}(x, t_{1}) = s = i^{L}(x, t_{2}), i^{L}(x, t_{1}') = s' = i^{L}(x, t_{2}').
\]

The set  $U$, as defined above,  is an open set in $R_{L}$ containing $(x_{1}, x_{2})$. To see this
set satisfies the necessary conditions, we let $x'$ be in $U_{1}$ and suppose that 
$\varphi^{b(x', t_{1}, t_{1}', t_{1}'')}(x')$ is in $B^{2}$. From the conditions, we have 
$b(x, t_{1}, t_{1}', t_{1}'') =  b(x', t_{1}, t_{1}', t_{1}'')$ and, for all
$ l \geq 1$,  $i^{l}(x', t_{1}) \neq i^{l}(x', t_{1}')$. On the other hand, we also have
\begin{eqnarray*}
i^{L}(x', t_{1}) = i^{L}(x, t_{1}) & = s =  & i^{L}(x, t_{2}) = i^{L}(x', t_{2}), \\
i^{L}(x', t_{1}') = i^{L}(x, t_{1}') & = s' =  & i^{L}(x, t_{2}') = i^{L}(x', t_{2}')
\end{eqnarray*}
From which it follows that
\begin{eqnarray*}
i^{l}(x', t_{1}) & = & i^{l}(x',t_{2}), \\
i^{l}(x', t_{1}') & = & i^{l}(x',t_{2}'),
\end{eqnarray*}
for all $l \geq L$. This then implies that 
\[
i^{l}(x', t_{2}) \neq  i^{l}(x',t_{2}'),
\]
for all $l \geq L$. From Lemma \ref{lemma:6_B}, we have $\varphi^{b(x', t_{2}, t_{2}', t_{2}'')}(x')$ is also in $B^{2}$.
\end{proof}

\begin{prop}
\label{prop:6_K2_transverse}
$K^{2}$ is transverse to $R \mid B^{2}$.
\end{prop} 

\begin{proof}
From part 4 of Lemma \ref{lemma:6_B} and Lemma \ref{lemma:6_K_compact},  $K^{2}$ is a compact, \'{e}tale equivalence relation on $B^{2}$. Next,
for any  $\varphi^{b(x,t, t', t'')}(x)$ in $B^{2}$,  by part 2 of Lemma \ref{lemma:6_B}
we have $i^{l}(x,t) \neq i^{l}(x, t')$, for every $ l \geq 1$. 
It follows from part 2 of Lemma \ref{lemma:6_i3} that $i_{3}(x, t, t', t'') = (i(x,t), i(x, t'), t_{1})$ and also that \linebreak 
$i_{3}(x, t', t, t'') = (i(x,t'), i(x, t), t_{1}')$. for some $t_{1}, t_{1}'$.  It then follows that \linebreak
$b(x,t, t', t'')$ is in $\Z^{2}_{1}(x, i(x,t), i(x,t'), t_{1})$ while 
$b(x,t', t, t'')$ is in \linebreak 
$\Z^{2}_{1}(x, i(x,t'), i(x,t), t_{1}')$. Continuing inductively, we see that, for every $l \geq 1$, there are $t_{l}, t_{l}'$ such that $b(x,t, t', t'')$ is in  
$\Z^{2}_{l}(x, i^{l}(x,t), i^{l}(x,t'), t_{l})$ while $b(x,t', t, t'')$ is in $\Z^{2}_{l}(x, i^{l}(x,t'), i^{l}(x,t), t_{l}')$. From this, we see that  \linebreak 
$(\varphi^{b(x,t, t', t'')}(x), \varphi^{b(x,t', t, t'')}(x))$ is not in $R_{l}$, for any $l \geq 1$. Thus $K^{2}$ has trivial intersection with $R$.

Moreover, if we have another point $\varphi^{b(x,s, s', s'')}(x)$ in $B^{2}$ such that
\linebreak 
$(\varphi^{b(x,s, s', s'')}(x), \varphi^{b(x,t, t', t'')}(x))$ is in $R$, then for some $l \geq 1$,
 $b(x,s, s', s'')$ and $b(x, t, t', t'')$ are in the same set $\Z^{2}_{l}(x, u, u', u'')$. From the argument above, we see that $u=i^{l}(x,t) = i^{l}(x,s), u'= i^{l}(x,t') =i^{l}(x,s')$ and then $b(x,t',t,t'')$ is in $\Z^{2}_{l}(x, u', u, u_{1}'')$ and 
$b(x,s',s,s'')$ is in $\Z^{2}_{l}(x, u', u, u_{2}'')$.
It follows from part 3 of  Lemma \ref{lemma:6_i3}, that $b(x,s', s, s'')$ and $b(x, t', t, t'')$ are in the same set   $\Z^{2}_{l+1}(x, v, v', v'')$. This shows that the map 
\begin{eqnarray*}
  &  ( (\varphi^{b(x,s, s', s'')}(x), \varphi^{b(x,t, t', t'')}(x)), (\varphi^{b(x,t, t', t'')}(x), \varphi^{b(x,t', t, t'')}(x) ) )
 &      \\ 
\mapsto & (  ( \varphi^{b(x,s, s', s'')}(x), \varphi^{b(x,s', s, s'')}(x)),
 (\varphi^{b(x,s', s, s'')}(x),\varphi^{b(x,t', t, t'')}(x) ) )&
\end{eqnarray*}
is a bijection between  $(R | B^{2}) \times_{B^{2}} K^{2}$ and $K^{2} \times_{B^{2}} (R | B^{2})$, as desired. It is easy to verify that it is a homeomorphism.
This completes the proof.
\end{proof}

\begin{prop}
\label{prop:6_B3_etale}
$B^{3}$ is an \'{e}tale subset of $B^{2}$ for the relation $(R \mid B^{2} ) \vee K^{2}$.
\end{prop} 

\begin{proof}
First, we note that the set 
\[
(R \mid B^{2} ) \vee K^{2} \setminus (R \mid B^{2})
\]
consists of those pairs 
 $(\varphi^{b(x,s, s', s'')}(x), \varphi^{b(x,t', t, t'')}(x))$, where $x$ is in $X$,
$(s, s', s'')$ and $ (t, t', t'')$ are in $\mathcal{T}_{0}^{3}(x)$ such that 
$\varphi^{b(x,s, s', s'')}(x), \varphi^{b(x,t, t', t'')}(x)$ are in $B^{2}$ and for some $L \geq 0$, 
$b(x,s, s', s'')$ and 
$b(x,t, t', t'')$ are in  $\Z^{2}_{L}(x, u, u', u'')$, for some $(u, u', u'')$ in 
$\mathcal{T}_{L}^{3}(x)$.

If we also consider the case that $\varphi^{b(x,s, s', s'')}(x)$ and $\varphi^{b(x,t', t, t'')}(x)$ 
are in $B^{3}$,
it follows from Lemma \ref{lemma:6_B}, that $i^{l}(x,s), i^{l}(x,s'), i^{l}(x,s'')$ are all distinct and 
$b(x, s, s', s'')$ is in $\Z^{2}_{l}(x, i^{l}(x,s), i^{l}(x,s'), i^{l}(x,s''))$, and also that
$i^{l}(x,t'), i^{l}(x,t)$, 
and $ i^{l}(x,t'')$ are all distinct and 
$b(x, t', t, t'')$ is in \linebreak 
$\Z^{2}_{l}(x, i^{l}(x,t'), i^{l}(x,t), i^{l}(x,t''))$, for all $l \geq 1$.
It follows that
\[
i^{L}(x, s) = i^{L}(x, t), i^{L}(x, s') = i^{L}(x, t'), i^{L}(x, s'') = i^{L}(x, t''),
\]
and hence,
\[
i^{l}(x, s) =   i^{l}(x, t), 
i^{l}(x, s')  =  i^{l}(x, t'), 
i^{l}(x, s'')  =  i^{l}(x, t''),
\]
for all $l \geq L$.

Now, using the regularity of $\mathcal{T}_{l}$, $0 \leq l \leq L$, we choose a clopen neighbourhood  $U_{0}$
of $x$ such that $(s, s', s''), (t, t', t'')$ are in $\mathcal{T}_{0}^{3}(x')$, for every $x'$
in $U_{0}$, and
\begin{eqnarray*}
b(x, s, s', s'') = b(x', s, s', s''), & b(x, s', s, s'') = b(x', s', s, s''), &\\
b(x, t, t', t'') = b(x', t, t', t''), & b(x, t', t, t'') = b(x', t', t, t''), & \\
i^{l}(x, s) = i^{l}(x', s), &  i^{l}(x, s') = i^{l}(x', s'), &  \\
i^{l}(x, s'') = i^{l}(x', s''), & i^{l}(x, t) = i^{l}(x', t), &  \\ 
i^{l}(x, t') = i^{l}(x', t'), & i^{l}(x, t'') = i^{l}(x', t''), & \\
\end{eqnarray*}
for all $0 \leq l \leq L$ and all $x'$ in $U_{0}$.

We consider 
\[
U = \{ (\varphi^{b(x,s,s',s'')}(x'), \varphi^{b(x,t',t,t'')}(x')) \mid x' \in U_{0} \cap B^{2} \},
\]
which is a neighbourhood of our point in  $(R \mid B^{2} ) \vee K^{2}$.
If we choose a point in $U$ with $\varphi^{b(x,s,s',s'')}(x')$ in $B^{3}$, it follows
that $i^{l}(x',s), i^{l}(x',s'), i^{l}(x',s'')$ are all distinct, for all $l \geq 0$.
We have
\[
i^{L}(x',s) = i^{L}(x,s) = i^{L}(x,t)  = i^{L}(x',t) ,
\]
from which it follows that
\[
i^{l}(x',s) = i^{l}(x',t),
\]
for all $l \geq L$. Similarly, we have 
\[
i^{l}(x',s') = i^{l}(x',t'), i^{l}(x',s'') = i^{l}(x',t''),
\]
for all $l \geq L$ and it follows that the second point, $\varphi^{b(x,t',t,t'')}(x')$ is also
in $B^{3}$. 

A similar proof holds for pairs in $R | B^{2}$. This completes the proof.
\end{proof}\

\begin{prop}
\label{prop:6_K3_transverse}
$K^{3}$ is transverse to $((R \mid B^{2} ) \vee K^{2}) \mid B^{3}$.
\end{prop} 

\begin{proof}
A typical
element of $(R \mid B^{2} ) \vee K^{2}$ is in one of two forms:
\[
(\varphi^{b(x,t,t',t'')}(x), \varphi^{b(x,s,s',s'')}(x)), (\varphi^{b(x,t,t',t'')}(x), \varphi^{b(x,s',s,s'')}(x))
\]
for some $x$ in $X$, $(t,t',t''), (s,s',s'')$ in $\mathcal{T}_{0}^{3}(x)$ such that \linebreak
$b(x, t, t', t''), b(x, s, s', s'')$ are in the same set $\Z^{2}_{L}(x, u, u', u'')$,
for some $L \geq 0$ and $(u, u', u'')$ in $\mathcal{T}_{L}^{3}(x)$. We consider
only the second case, the first is similar. If, in addition,
we consider points which are in $B^{3}$, then we have
\[
i^{L}(x, t) = u = i^{L}(x, s), i^{L}(x, t') = u' = i^{L}(x, s'), i^{L}(x, t'') = u'' = i^{L}(x, s'').
\]
First, notice that for such a pair, the $(s',s,s'')$ cannot be an even permutation of
$(t, t', t'')$ and hence, $K^{3}$ has trivial intersection with 
$(R \mid B^{2} ) \vee K^{2}$.
Secondly, if we take the pairs $(\varphi^{b(x,t'',t,t')}(x), \varphi^{b(x,t,t',t'')}(x))$, 
which is in $K^{3}$,  and $(\varphi^{b(x,t,t',t'')}(x), \varphi^{b(x,s',s,s'')}(x))$, as above,
it follows that \linebreak
$ (\varphi^{b(x,t'',t,t')}(x), \varphi^{b(x,s,s'',s')}(x))$ is in
$(R \mid B^{2} ) \vee K^{2}$, while \linebreak
$(\varphi^{b(x,s,s'',s')}(x), \varphi^{b(x,s',s,s'')}(x))$ is in $K^{3}$. This establishes the bijection between
$K^{3} \times_{B^{3}} (((R \mid B^{2}) \vee K^{2}) | B^{3})$ and 
$(((R \mid B^{2}) \vee K^{2}) | B^{3}) \times_{B^{3}} K^{3}$, at least in the second case above. It is clear that this is
a homeomorphism. 
The details for the first case are left to the reader.
\end{proof}

What remains is for us to show that the orbit relation, $R_{\varphi}$, is 
generated, as an equivalence relation, by $R$, $K^{2}$ and $K^{3}$.
We begin with a lemma.

\begin{lemma}
\label{lemma:6_R_class}
Let $x$ be in $X$.
\begin{enumerate}
\item
Let $t$ be in $\mathcal{T}_{l}(x)$. The elements  $\varphi^{n}(x), n \in \Z^{2}  \cap'  t $ are all in the same
equivalence class in $R_{l+1}$, unless there exists $(t, t', t'')$ in
$\mathcal{T}_{l}^{3}(x)$ such that $i(x,t), i(x,t'), i(x,t'')$ are all distinct.
In this case, $ \Z^{2} \cap' t$ is contained in the union of
$\Z^{2}_{l+1}(x, i(x,t), i(x,t'), i(x,t''))$ and \linebreak
 $\Z^{2}_{l+1}(x, i(x,t), i(x,t''), i(x,t'))$.
\item
Let $(t, t', t'')$ be in $\mathcal{T}_{l}^{3}$(x). For 
$m \in \Z^{2}_{l}(x, t, t', t'')$ and $ n \in \Z^{2}_{l}(x, t, t'', t')$, $\varphi^{n}(x), \varphi^{m}(x)$
 are in the same 
$R_{l+1}$-equivalence class unless \linebreak 
$i(x,t), i(x,t'), i(x,t'')$ are all distinct.
\end{enumerate}
\end{lemma}

\begin{proof}
The first fact follows from a  close examination of the definition of $i_{3}$, noticing that
the result depends on $t', t''$ only if $i(x,N(x,t))$ contains three
elements of $\mathcal{T}_{l+1}(x)$ and the rest of the conclusion follows easily.

The second fact also follows from a  close examination of the definition of $i_{3}$, 
noticing that the result only depends on the order of $t', t''$ in the case $i(x,t), i(x,t')$ and $i(x,t'')$ are all distinct.
Again the conclusion is clear.
\end{proof}

\begin{prop}
\label{prop:6_generate}
We have 
\[
R \vee K^{2} \vee K^{3} = R_{\varphi}.
\]
\end{prop}

\begin{proof}
Let $x$ be in $X$. We will show that the result holds for the $R_{\varphi}$-equivalence class of $x$.
For a fixed $r > 0$, we define
\[
n_{r} = \lim_{l \rightarrow \infty} \# \{ t \in \mathcal{T}_{l}(x) \mid t \cap B(0,r) \neq \emptyset \}.
\]
First observe that, for any fixed $l \geq 0$, the map $i(x, \cdot)$ is a surjection from 
the set in the expression above for $l$, to that for $l+1$, hence the sequence above
is a decreasing sequence of positive integers and, in particular, the limit exists.
Secondly, for $l > 2r + 1$, by part 2 of Theorem \ref{thm:5_refine},  $\mathcal{T}_{l}(x)$ is $2r$-separated.
Hence if two elements of $\mathcal{T}_{l}(x)$ meet $B(0,r)$, then they meet each other. Since $\mathcal{T}_{l}(x)$ is
simplicial, it follows that at most three elements meet $B(0,r)$, for large values of $l$, 
and 
hence, $n_{r} \leq 3$. Finally, notice that if $r' > r$, then $n_{r'} \geq n_{r}$.

We first consider the case that, for every $r > 0$, $n_{r} =1$. Given any $n$ in $\Z^{2}$,
using the fact that $n_{r}=1$, 
we may find $l \geq 0$ such that there is a unique $t$ in $\mathcal{T}_{l}(x)$ which meets
$B(0, |n| +1)$. From this we conclude that $0, n \in \Z^{2} \cap' t $. Now suppose
that $x$ and $\varphi^{n}(x)$ are not in the same $R_{l+1}$-class. From the first part of \ref{lemma:6_R_class}, 
we find $(t,t',t'')$ in $\mathcal{T}_{l}^{3}(x)$ such that $i(x,t), i(x,t'),i(x, t'')$
are all distinct and $0 \in \Z^{2}_{l+1}(x, i(x,t),i(x,t'),i(x, t''))$, 
$n \in \Z^{2}_{l+1}(x, i(x,t),i(x,t''),i(x, t'))$. If $\varphi^{0}(x), \varphi^{n}(x)$ 
are not in the same 
$R_{l+k}$ class for any $k \geq 1$, then by the second part of \ref{lemma:6_R_class}, we have 
 $i^{k}(x,t), i^{k}(x,t')$ and $i^{k}(x, t'')$ are all distinct, for $k \geq 1$.
But then, we may find $r'>r$ such that $B(0,r')$ meets each of $t, t', t''$.
And we have $B(0,r')$ meets each of $i^{k}(x,t), i^{k}(x,t')$ and $i^{k}(x, t'')$, for every $k \geq 1$,
It follows then that $n_{r'} =3$. As this is impossible, we conclude that 
$x, \varphi^{n}(x)$ are in the same $R_{l+k}$ class for some $k \geq 0$. Since $n$ was arbitrary, we conclude
that the $R$ class of $x$ is the same as its $R_{\varphi}$-class.

Secondly, we consider the case that, for some $r$, $n_{r}=2$, but for all $r'$, $n_{r'} \neq 3$.
Let $L$ be such that $\mathcal{T}_{L}(x)$ contains two elements, $t$ and $t'$ which meet
$B(0,r)$. As $n_{r} = 2$, $i^{L+k}(x, t), i^{L+k}(x, t')$ are distinct for all $k \geq 1$.
Since the sets $t$ and $t'$ are closed and $B(0,r)$ is connected, $t \cap t'$ is non-empty.
If we write $t$ and $t'$ as unions of elements of $\mathcal{T}_{0}(x)$, and using the distributive
law for intersections and unions, we find distinct elements $s, s'$ of $\mathcal{T}_{0}(x)$ such that
$s \cap s'$ is non-empty, $i^{L}(x,s) =t, i^{L}(x,s') =t'$. This means that
$i^{l}(x,s) \neq i^{l}(x, s')$, for all $l \geq 1$. Find $s''$ in $\mathcal{T}_{0}(x)$
such that $(s, s', s'')$ is in $\mathcal{T}_{0}^{3}(x)$. The pair
$(\varphi^{b(x,s,s',s'')}(x), \varphi^{b(x,s',s,s'')}(x))$ is in $K^{2}$. Now, let $n$ be
in $\Z^{2}$. We will show that either $(\varphi^{b(x,s,s',s'')}(x), \varphi^{n}(x))$ or
$(\varphi^{b(x,s',s,s'')}(x), \varphi^{n}(x))$ is in $R$ and this will complete the proof.
Let $r' = \max \{ r, |n|+1 \}$. For sufficiently large $k \geq 1$, $B(0,r')$ meets only 
$i^{k}(x, t) = i^{L+k}(x, s)$ and $i^{k}(x, t') = i^{L+k}(x, s')$. Let us suppose
that $n \in' i^{k}(x, t) = i^{L+k}(x, s)$. It follows from exactly the same argument as in the first case, that if
$\varphi^{n}(x)$ and $\varphi^{b(x, s, s', s'')}(x)$ are  not in the same 
$R$ equivalence class, then we may find $r''$ such that $n_{r''} = 3$, which
is not possible. This completes the proof in this case.

The final case to consider is that $n_{r} = 3$, for some $r > 0$. Find $L \geq 1$ such that
$\mathcal{T}_{L}(x)$ has exactly three elements, say $t, t', t''$ which meet $B(0,r)$.
Since $n_{r}=3$, $i^{k}(x, t), i^{k}(x, t'), i^{k}(x, t'')$ are all distinct, for 
all $k \geq 1$. For $L > 2r+1$, $\mathcal{T}_{L}(x)$ is $2r$-separated. This implies that each of the intersections $t \cap t'$, $t \cap t''$ and $t' \cap t''$ is non-empty, which, in turn, means that $t \cap t' \cap t''$ must be non-empty, since $\mathcal{T}_{L}(x)$ is simplicial. Next, we write each of $t, t', t''$ as a union 
of elements of $\mathcal{T}_{0}(x)$ and, using the distributive law for intersections
and unions, we may find $(s, s', s'')$ in $\mathcal{T}_{0}^{3}(x)$ with
$i^{L}(x, s) = t, i^{L}(x, s') = t', i^{L}(x, s'') = t''$. If $\sigma$ is any permutation of
$\{ s, s', s'' \}$, the pair 
$(\varphi^{b(x,s, s', s'')}(x), \varphi^{b(x,\sigma(s), \sigma(s'), \sigma(s''))}(x))$ is in $K^{2} \vee K^{3}$.
It remains for us to show that, for any $n$ in $\Z^{2}$, there is a $\sigma$ such that
$(\varphi^{n}(x), \varphi^{b(x,\sigma(s), \sigma(s'), \sigma(s''))}(x))$ is in $R$. We let 
$r' = \max \{ |n|+1, r \}$. Since $n_{r'} = 3$, we may find $k \geq 1$ such that 
$B(0, r')$ meets only three elements of $\mathcal{T}_{L+k}(x)$ and these must be
$i^{L+k}(x, s) = i^{k}(x, t), i^{L+k}(x, s') = i^{k}(x, t'), i^{L+k}(x, s'') = i^{k}(x, t'')$.
Then $n$ must be in one of them, say in $i^{L+k}(x, s) = i^{k}(x, t)$. Notice that 
$i^{L+k}(x, s')$ and $i^{L+k}(x, s'')$ are both in $N(x, i^{L+k}(x, s))$, and so
from the definition of $i_{3}$, we see that $n$ is in either
$\Z^{2}_{L+k+1}(x, i^{L+k+1}(x, s), i^{L+k+1}(x, s'), i^{L+k+1}(x, s''))$ or in 
\linebreak 
$\Z^{2}_{L+k+1}(x, i^{L+k+1}(x, s), i^{L+k+1}(x, s''), i^{L+k+1}(x, s'))$. In the former
case, $n$ and $b(x,s,s', s'')$ are both in $\Z^{2}_{L+k+1}(x, i^{L+k+1}(x, s), i^{L+k+1}(x, s'), i^{L+k+1}(x, s''))$
and hence $(\varphi^{n}(x), \varphi^{b(x,s,s',s'')}(x))$ is in 
$R_{L+k+1}$. The other case is similar and we have $(\varphi^{n}(x), \varphi^{b(x,s,s'',s')}(x))$ is in 
$R_{L+k+1}$. This completes the last case.
\end{proof}

We need to show  all the hypotheses of the absorption theorem are satisfied before proving the main result.

\begin{prop}
\label{prop:6_R_minimal}
The AF-relation $R$ is minimal.
\end{prop}

\begin{proof}
We first claim that for any $x$ in $X$ and $l \geq 0$, the set 
\[
\{ n \in \Z^{2} \mid (x, \varphi^{n}(x)) \in R_{l+1} \}
\]
contains a ball of radius $l$. The point $0$ is in $\Z^{2} \cap' t$ for some
$t$ in $\mathcal{T}_{l}(x)$. If there is  no $(t, t', t'')$ in 
$\mathcal{T}_{l}^{3}(x)$ with $i(x,t), i(x,t'), i(x,t'')$ are  all
distinct, it follows from \ref{lemma:6_R_class} that 
\[
\{ (x, \varphi^{n}(x)) \mid n \in 
\Z^{2} \cap' t \}
\]
is contained in $R_{l+1}$ and the conclusion follows since $\mathcal{T}_{l}(x)$ has
capacity $l$. By Theorem \ref{thm:5_refine}, in the case that there exists $(t, t', t'')$ in 
$\mathcal{T}_{l}^{3}(x)$ such that $i(x,t), i(x,t'), i(x,t'')$ are  all
distinct, see again from \ref{lemma:6_R_class} that $0$ is in 
either $\Z^{2}_{l+1}(x, i(x, t), i(x, t'), i(x, t''))$ or
 $\Z^{2}_{l+1}(x, i(x, t), i(x, t''), i(x, t'))$. Let us suppose the former. 
We appeal to condition 6 of \ref{thm:5_refine}, using 
$t_{1}= i(x,t), t_{2} = i(x,t'), t_{3} = i(x, t'')$, we find $s$ in $\mathcal{N}^{2}(x, t_{1}, t_{2}, t_{3})$. Condition 6 of  \ref{thm:5_refine} and Lemma \ref{lemma:6_R_class} imply  that
$\Z^{2} \cap' s$ is contained in  $\Z^{2}_{l+1}(x,t_{1}, t_{2}, t_{3})$. 
The conclusion again follows since $\mathcal{T}_{l}(x)$ has
capacity $l$.

We wish to show that the $R$-equivalence class of $x$ is dense in $X$. 
For each $l \geq 0$, we find $n_{l}$ in $\Z^{2}$ such that 
$(x,\varphi^{n}(x))$ is in $R_{l+1}$, for all $n$ in $B(n_{l}, l)$. 
By passing to a subsequence, we may assume that $\varphi^{n_{l}}(x)$
converges to a point, say $x'$, in $X$. Now let $n$ be in $\Z^{2}$. 
The sequence $\varphi^{n_{l}+n}(x)$ converges to $\varphi^{n}(x')$. 
For $ l > |n|$,  $n_{l}+n$ is in $B(n_{l}, l)$ and so
$(x, \varphi^{n_{l}+n}(x))$ is in $R$. It follows that $\varphi^{n}(x')$
is in the closure of the $R$-equivalence class of $x$. But since $n$ was 
arbitrary, the entire $\varphi$-orbit of $x'$ is in the closure
of the $R$-equivalence class of $x$. Since the $\varphi$-orbit of any point 
is dense, we conclude that the $R$-equivalence class of $x$ is dense in $X$.
\end{proof}

\begin{prop}
\label{prop:6_B2_measure}
For any $R$-invariant probability measure $\mu$ on $X$, we have
\[
\mu(B^{2}) = \mu(B^{3}) = 0.
\]
\end{prop}

\begin{proof}
We consider $B^{2}$ first, the conclusion for $B^{3}$ follows since it is a subset
of $B^{2}$. We consider
$(t_{1}, t_{2}, t_{3})$ and $(t_{1}, t_{2}, t_{3}' )$ in 
$\mathcal{T}_{l+1}^{3}(x)$ with $t_{3} \neq t_{3}'$. We will show that
there is a collection of functions, $ 1 \leq i \leq 2^{l-2}$, 
\[
\eta_{i}: B^{2}_{l+1}(x, t_{1}, t_{2}) \cap \Z^{2}_{l+1}(x, t_{1}, t_{2}, t_{3})
\rightarrow \Z^{2}_{l+1}(x, t_{1}, t_{2}, t_{3})
\]
with pairwise disjoint ranges, and for each $i$,
\[
\{ (\varphi^{n}(x), \varphi^{\eta_{i}(n)}(x)) \mid x \in X, n \in B^{2}_{l+1}(x, t_{1}, t_{2}) \} \subset R_{l+1}.
\]
 Moreover, these may be chosen
in a $\varphi$-regular way (although Borel is sufficient). It follows from this, that
if $\mu$ is any $R_{l+1}$-invariant probability measure on $X$, we have 
\[
\mu(B^{2}_{l+1}) \leq 2^{-l+2},
\]
and the conclusion follows at once.

To construct the functions, we proceed as follows. Since $K^{2}_{0}$ is compact, 
there is $L \geq 1$ such that, for any $x$ in $X$,
and $(s,s',s'')$ in $\mathcal{T}_{0}^{3}(x)$, we have 
 $d(b(x, s, s', s''), b(x, s', s, s'')) \leq L $. (Since $b(x,s,s',s'')$ and and $b(x,s',s,s'')$ lie in $s$ and $s'$, respectively, and these sets intersect, the distance between them is at most $2diam(\mathcal{T}_{0})$.) Consider $ l \geq L+1$. If $b$
is a point in $B^{2}_{l}(x, t, t')$, for some $(t, t')$ in $\mathcal{T}_{l}^{2}(x)$, 
then there is a $b'$  in $B^{2}_{l}(x, t', t)$ with $d(b, b') \leq l - 1$. It follows
that 
\[
B^{2}_{l}(x,t,t') \subset \{ n \in \Z^{2} \cap' t  \mid d( n, \R^{2} \setminus t)  \leq l - 1 \}.
\]
The domain is contained in the union of sets
$B^{2}_{l}(x, t, t') \cap \Z^{2}_{l}(x, t, t', t'')$, where $(t, t', t'')$ are in 
 $\mathcal{T}_{l}^{3}(x)$ with $i_{3}(x, t, t', t'') = (t_{1}, t_{2}, t_{3}).$ 
For our domain to have non-empty intersection with such a set requires
$i(x, t) = t_{1}, i(x, t') = t_{2}$. First consider the case $i(x, t'') = t_{3}$.
From condition 6 of \ref{thm:5_refine}, each such triples $(t, t', t'')$
is in $\mathcal{N}^{3}(x, t_{1}, t_{2}, t_{3})$. The number of such triple is
less than $\mathcal{N}^{2}(x, t_{1}, t_{2}, t_{3})$.
So we may find, for each $(t, t', t'')$, $s$ in $\mathcal{N}^{2}(x, t_{1}, t_{2}, t_{3})$
such that this function (which we do not name) is injective. 
With the observation above in the first paragraph and condition 7 of
\ref{thm:5_refine}, we know that
\[
2^{l-1} \# (B^{2}_{l}(x, t, t') \cap \Z^{2}_{l}(x, t, t', t'')) \leq  \#( \Z^{2} \cap' s ).
\]
This means we may find $2^{l-2}$ injective functions with pairwise disjoint ranges
from the first set into the second. Moreover, the complement of
all their ranges still contains at least $\#( \Z^{2} \cap' s)/2$ points.
For  $s$ in $\mathcal{N}^{2}(x, t_{1}, t_{2}, t_{3})$, we have $d(s,t_{3}) < d(s, t_{3}')$, by definition. Then for any $s', s''$ with $(s,s',s'')$ in $\mathcal{T}_{l}^{3}(x)$, $i_{3}(x,s,s',s'') = (t_{1}, t_{2}, t_{3})$.
From this, it follows that  $ \Z^{2} \cap' s $
is contained in $\Z^{2}_{l+1}(x, t_{1}, t_{2}, t_{3})$.

We may write $B^{2}_{l+1}(x, t_{1}, t_{2}) \cap \Z^{2}_{l+1}(x, t_{1}, t_{2}, t_{3})$ as a union of two sets. The first is the union of all sets 
$B^{2}_{l+1}(x, t_{1}, t_{2}) \cap \Z^{2}_{l+1}(x, t, t', t'')$, over all $(t, t', t'')$ such that $i(x,t) = t_{1}, i(x,t') = t_{2}, i(x,t'') = t_{3}$. We have already dealt with this part in the proof. The second set is the union of sets of the form
$B^{2}_{l+1}(x, t_{1}, t_{2}) \cap (\Z^{2} \cap' t)$, over all $t$ such that $i(x,t) = t_{1}$, $i(x, N(x,t)) = \{ t_{1}, t_{2} \}$. For any one of these sets, we observe 
\[
 B^{2}_{l+1}(x, t_{1}, t_{2}) \cap (\Z^{2} \cap' t) \subset \{ n \in \Z^{2} \cap' t \mid d(n, \R^{2} \setminus t) \leq l - 1 \},
\]
and so by condition 7 of Theorem \ref{thm:5_refine}, 
\[
 2^{l-1} \# (B^{2}_{l+1}(x, t_{1}, t_{2}) \cap (\Z^{2} \cap' t)) \leq \#  (\Z^{2} \cap' t).
\]
Hence we may define $2^{l-2}$ injective functions from $B^{2}_{l+1}(x, t_{1}, t_{2}) \cap (\Z^{2} \cap' t)$ to $ \Z^{2} \cap' t$ with pairwise disjoint images. Moreover, these images may be chosen to be disjoint from those obtained in the first case.
It is clear that all of these functions may be locally derived in the appropriate sense. This completes the proof.
\end{proof}

At this point, we have completed the proof of Theorem \ref{thm:affapp}.


\end{document}